\documentclass{article}[12pt]
\input epsf.tex
\usepackage{latexsym}
\usepackage{amsfonts,amssymb}
\usepackage{epsfig}
\usepackage{psfig}

\newtheorem{theorem}{Theorem}[section]
\newtheorem{lemma}[theorem]{Lemma}
\newtheorem{corollary}[theorem]{Corollary}

\newtheorem{question}[theorem]{Question}

\newtheorem{proposition}[theorem]{Proposition}

\newtheorem{convention}[theorem]{Convention}

\def\Mod{\mbox{\rm{Mod}}}
\def\Homeo{\mbox{\rm{Homeo}}}

\def\Stab{\mbox{\rm{Stab}}}

\def\PSL{\mbox{\rm{PSL}}}
\def\SO{\mbox{\rm{SO}}}
\def\SU{\mbox{\rm{SU}}}
\def\Aff{\mbox{\rm{Aff}}}

\def\Aut{\mbox{\rm{Aut}}}
\def\prj{{\mathbb P}}
\def\ML{{\cal ML}}

\def\Im{\mbox{\rm{Im}}}
\def\vcd{\mbox{\rm{vcd}}}

\def\intr{\mbox{\rm{int}}}

\begin{document}

\title{A combination theorem for Veech subgroups of the mapping class
group}
\author{Christopher J. Leininger \thanks{This work was partially
supported by an N.S.F. postdoctoral fellowship} \& Alan W. Reid 
\thanks{This work was partially supported by an N.S.F. grant}}

\maketitle

\begin{abstract}
In this paper we prove a combination theorem for Veech subgroups of the
mapping class group analogous to the first Klein-Maskit combination
theorem for Kleinian groups in which two Fuchsian subgroups are
amalgamated along a parabolic subgroup.
As a corollary, we construct subgroups of the mapping class group (for all
genus at least $2$), which are isomorphic to non-abelian closed surface
groups in which all but one conjugacy class (up to powers) is
pseudo-Anosov.
\end{abstract}

%%%%%%%%%%%%%%%%%%%%%%%%%%%%%%%%%%%%%%%%%%%
%%%%%%%%%%%%%%%%%%%%%%%%%%%%%%%%%%%%%%%%%%%
%%%%%%%%%%%%%%%%%%%%%%%%%%%%%%%%%%%%%%%%%%%
\section{Introduction} \label{introsect}
%%%%%%%%%%%%%%%%%%%%%%%%%%%%%%%%%%%%%%%%%%%
%%%%%%%%%%%%%%%%%%%%%%%%%%%%%%%%%%%%%%%%%%%
%%%%%%%%%%%%%%%%%%%%%%%%%%%%%%%%%%%%%%%%%%%

For $R$ a compact oriented surface, possibly with boundary,
the {\em mapping class group} of $R$ is defined to be the group of isotopy
classes of orientation preserving homeomorphisms
$$\Mod(R) =  \Homeo^{+}(R) / \Homeo_{0}(R),$$
where, $\Homeo^{+}(R)$ is the group of orientation preserving
homeomorphisms of $R$ and $\Homeo_{0}(R)$ are those homeomorphisms
isotopic to the identity.
We will write $\phi$ to denote a homeomorphism and its isotopy class, when
no confusion can arise, or when the distinction is unimportant.
An isotopy class of a homeomorphism is called an {\em automorphism} of
$R$. We specify no boundary behavior when $R$ has boundary.

There have been many analogies made between $\Mod(R)$ 
and lattices in Lie groups.  In particular, the study of the
subgroup structure of $\Mod(R)$, both finite and infinite index, 
has many parallels in the theory of lattices. For example, the question
of whether Property T holds for $\Mod(R)$ ($R$ a closed orientable surface
of genus $\geq 3$) and whether $\Mod(R)$ has a version of the Congruence
Subgroup Property, or towards the other extreme, whether there are finite
index subgroups of $\Mod(R)$ that surject onto ${\mathbb Z}$ are questions that have received much attention recently (see for example \cite{Ha}, \cite{Iv}, and \cite{FLM} for more on this).

The point of view of this paper is similar to \cite{MP} and \cite{FMo} and is motivated by results and 
constructions in 3-manifold topology and Kleinian groups.
A powerful tool in the theory of Kleinian groups are the
Klein-Maskit combination theorems.
Our purpose here is to take the first step in the development of 
combination theorems for subgroups of $\Mod(R)$.

Before we state our main theorem, we describe the motivating example from
Kleinian groups.
Suppose $G_{1},G_{2} < \PSL_{2}({\mathbb R})$ are two finitely generated
Fuchsian groups of finite co-area, and
$G_{0}$ is a maximal parabolic subgroup of both.
If $h \in \PSL_{2}({\mathbb C})$ is another parabolic centralizing
$G_{0}$, then we may form the amalgamated product
$$G = G_{1} *_{G_{0}} h G_{2} h^{-1}$$

The First Klein-Maskit Combination theorem tells us that if $h$ is
"sufficiently complicated" (see \S \ref{kleinsect}), then the natural
homomorphism of $G$ into $\PSL_{2}({\mathbb C})$ is injective and
every element is hyperbolic except those conjugate into an elliptic or
parabolic subgroup of either factor.  In the particular case that each
$G_{1}$ and $G_{2}$ are torsion free and for $i=1,2$, ${\Bbb H}^2/G_i$
have one cusp, $G$ is isomorphic to a non-abelian closed surface
subgroup of $\pi_{1}(M)$ with one accidental parabolic. This
has been very useful in the construction of essential closed surfaces
in cusped hyperbolic 3-manifolds (see \cite{CLR} and \cite{CL}).

Our main theorem states that we can do essentially the same thing in
$\Mod(R)$.
The Fuchsian subgroups in this setting are the {\em Veech subgroups} of
$\Mod(R)$.
A hyperbolic element of $\Mod(R)$ is a pseudo-Anosov automorphism.

A special case of our main Theorem \ref{combothrm} is given by the
following (see \S \ref{grpthrysect}-- \S \ref{mcgsect} for definitions).

\begin{theorem} \label{fakecombothrm}
Suppose $G(q_{1})$, $G(q_{2})$ are Veech subgroups of $\Mod(R)$ with
$G_{0}$ a maximal parabolic subgroup of each.
If $h \in \Mod(R)$ centralizes $G_{0}$ and is "sufficiently complicated",
then
$$G(q_{1}) *_{G_{0}} h G(q_{2}) h^{-1} \hookrightarrow \Mod(R)$$ is
injective.  Moreover, any element not conjugate into a
parabolic subgroup of either factor is pseudo-Anosov or a periodic element
conjugate into a factor.
\end{theorem}

One cusped Veech group lattices exist in most mapping class groups and in
particular we have

\begin{corollary} \label{sfcegrpcor}
For every closed surface $R$ of genus $g \geq 2$, there exist subgroups of
$\Mod(R)$ isomorphic to the fundamental group of a closed surface (of
genus $2g$) for which all but one conjugacy class of elements (up to
powers) is pseudo-Anosov.
\end{corollary}

That surface subgroups exist in most mapping class groups is already known
(see for example \cite{GH} for some discussion of this).
One new feature about our construction
is that ``most'' elements are pseudo-Anosov, unlike the other constructions.

These examples are of particular interest in connection with the following
question which has implications into the existence of negatively
curved $4$-manifolds that fiber over a surface (see the discussion at the 
end of \cite{GH} and Question 12.3 of \cite{Be}).

\begin{question}
Do there exist subgroups of the mapping class group isomorphic to the
fundamental group of a closed surface of genus at least $2$ for which all
non-identity elements are pseudo-Anosov?
\end{question}

The rest of the paper is organized as follows.
\S \ref{grpthrysect} gives the necessary definitions from combinatorial
group theory and ends with a more detailed description of the Kleinian
group example mentioned above.
In \S \ref{topcombsurfsect} we recall some of the tools used in studying
mapping class groups which we will need.
Veech groups are defined by a Euclidean cone structure on a surface which
we describe in \S \ref{qdsect}.
The first half of \S \ref{mcgsect} gives a brief recap of the
Nielsen-Thurston classification of surface automorphisms and provides the
theorems on pseudo-Anosov automorphisms we will need.
The end of this section describes Veech groups and provides examples.

Section \ref{combosect} contains the statement of Theorem \ref{combothrm}
(thus clarifying Theorem \ref{fakecombothrm}) and a proof of Corollary
\ref{sfcegrpcor} assuming this. The proof of Theorem \ref{combothrm} is
given in \S \ref{proofsect}, with results and notation needed in the proof
built up in \S \ref{sadtracsect} and \ref{hballsect}.
In \S \ref{sadtracsect} we exploit the structure of geodesic
representatives of simple closed curves in the Euclidean cone metrics
under consideration and use this to define sets which turn out to have
properties analogous to those used in the Kleinian group example.
Continuing to follow this example, we define in \S \ref{hballsect} those
sets which will play the part of horoballs, and
show that they do indeed have the necessary properties.
As well as proving Theorem \ref{combothrm} in \S \ref{proofsect}, 
we give some further generalizations of
Theorem \ref{combothrm} including a very  simple version of the Second
Klein-Maskit Combination theorem.
We end with some concluding remarks and questions in \S \ref{remarksect}.

%%%%%%%%%%%%%%%%%%%%%%%%%%%%%%%%%%%%%%%%%%
%%%%%%%%%%%%%%%%%%%%%%%%%%%%%%%%%%%%%%%%%%
%%%%%%%%%%%%%%%%%%%%%%%%%%%%%%%%%%%%%%%%%%
\section{Group theory} \label{grpthrysect}
%%%%%%%%%%%%%%%%%%%%%%%%%%%%%%%%%%%%%%%%%%
%%%%%%%%%%%%%%%%%%%%%%%%%%%%%%%%%%%%%%%%%%
%%%%%%%%%%%%%%%%%%%%%%%%%%%%%%%%%%%%%%%%%%

Since it will be useful in what follows we recall some basic facts about 
free products with amalgamations.

%%%%%%%%%%%%%%%%%%%%%%%%%%%%%%%%%%%%%%%%%%%%%%%%%%%%%%%%%%%%%
%%%%%%%%%%%%%%%%%%%%%%%%%%%%%%%%%%%%%%%%%%%%%%%%%%%%%%%%%%%%%
\subsection{Free product with amalgamation} \label{amalgsect}
%%%%%%%%%%%%%%%%%%%%%%%%%%%%%%%%%%%%%%%%%%%%%%%%%%%%%%%%%%%%%
%%%%%%%%%%%%%%%%%%%%%%%%%%%%%%%%%%%%%%%%%%%%%%%%%%%%%%%%%%%%%

Given groups $G_{0},G_{1},...,G_{P}$ and monomorphisms
$$\nu_{i}:G_{0} \rightarrow G_{i}$$
for each $i = 1,...,P$, the {\em free product of $\{ G_{i} \}_{i=1}^{P}$
amalgamated over $G_{0}$} is the group
\begin{equation} \label{amalgeqn}
G = G_{1} *_{G_{0}} \cdots *_{G_{0}} G_{P}
\end{equation}
defined by the universal property:
\begin{enumerate}
\item For each $i$, we have a monomorphism $\iota_{i}: G_{i} \rightarrow
G$ with $\iota_{i} \circ \nu_{i} = \iota_{i'} \circ \nu_{i'}$ for every
$i,i'$.
\item For any group $K$ and any collection of homomorphisms $\{ \eta_{i} :
G_{i} \rightarrow K \}$ with $\eta_{i} \circ \nu_{i} = \eta_{i'} \circ
\nu_{i'}$ for each $i,i'$, there exists a unique homomorphism $\eta:G
\rightarrow K$
such that $\eta \circ \iota_{i} = \eta_{i}$ for each $i$.
\end{enumerate}
This allows us to unambiguously identify $G_{0}$ and each of the $G_{i}$'s
as subgroups of $G$.

One can construct this group from the free product
$$G_{1} * \cdots * G_{P}$$
as the quotient by the normal subgroup generated by all elements of the
form $\nu_{i}(\phi) * \nu_{i'}(\phi^{-1})$ for $\phi \in G_{0}$, $i,i' \in
\{ 1,...,P \}$.
In particular, this gives rise to normal forms for the elements of $G$.
Namely, for every element $\phi \in G$ exactly one of the following holds:
\begin{enumerate}
\item $\phi$ is an element of $G_{0}$,
\item $\phi \in G_{i} \setminus G_{0}$ for some $i \in \{ 1,...,P \}$.
\item $\phi = \phi_{i_{1}} \cdots \phi_{i_{r}}$ where $\phi_{i_{j}} \in
G_{i_{j}} \setminus G_{0}$ with $i_{j} \in \{ 1,...,P \}$ and $i_{j} \neq
i_{j+1}$ for $j = 1,...,r-1$.
\end{enumerate}

%%%%%%%%%%%%%%%%%%%%%%%%%%%%%%%%%%%%%%%%%%%%%%%%%%%%%%%%%%
%%%%%%%%%%%%%%%%%%%%%%%%%%%%%%%%%%%%%%%%%%%%%%%%%%%%%%%%%%
\subsection{Proper interactive $P$-tuples} \label{pipsect}
%%%%%%%%%%%%%%%%%%%%%%%%%%%%%%%%%%%%%%%%%%%%%%%%%%%%%%%%%%
%%%%%%%%%%%%%%%%%%%%%%%%%%%%%%%%%%%%%%%%%%%%%%%%%%%%%%%%%%

Suppose we are now given subgroups $G_{0},G_{1},...,G_{P} < \Gamma$ with
$G_{i} \cap G_{i'} = G_{0} < \Gamma$, for every $i,i' \in \{ 1,...,P \}$
and $i \neq i'$.
Using the inclusions of $G_{0}$ into each $G_{i}$ we form the amalgamated
product $G$ as in (\ref{amalgeqn}).
The inclusions of each $G_{i}$ into $\Gamma$ and property 2 of the
amalgamated product implies the existence of a unique homomorphism $G
\rightarrow \Gamma$ extending these inclusions.

If $\Gamma$ acts on a set $X$, then we say that a $P$-tuple of subsets
$\Theta_{1},...,\Theta_{P} \subset X$ is a {\em proper interactive
$P$-tuple} for $G_{1},...,G_{P}$ if
\begin{enumerate}
\item $\Theta_{i} \neq \emptyset$ for each $i$,
\item $\Theta_{i} \cap \Theta_{i'} = \emptyset$ for each $i \neq i'$,
\item $G_{0}$ leaves $\Theta_{i}$ invariant for each $i$,
\item for every $\phi_{i} \in G_{i} \setminus G_{0}$ we have
$\phi_{i}(\Theta_{i'}) \subset \Theta_{i}$, for each $i' \neq i$, and
\item for every $i$, there exists $\theta_{i} \in \Theta_{i}$, such that
for every $\phi_{i} \in G_{i} \setminus G_{0}$, $\theta_{i} \not \in
\phi_{i}(\Theta_{i'})$ for any $i' \neq i$.
\end{enumerate}

The next proposition is proven using a standard ``ping-pong'' argument.
The case of two subgroups $G_{1},G_{2}$ amalgamated over $G_{0}$ is proved
in \cite{Mask}.

\begin{proposition} \label{pipprop}
Suppose $G_{0},G_{1},...,G_{P},\Gamma,X$ are as above and
$\Theta_{1},...,\Theta_{P} \subset X$ is a proper interactive $P$-tuple
for $G_{1},...,G_{P}$.
Then
$$G = G_{1} *_{G_{0}} G_{2} *_{G_{0}} \cdots *_{G_{0}} G_{P}
\hookrightarrow \Gamma$$
is an embedding.
\end{proposition}

\noindent
{\em Proof.} Given $\phi \in G \setminus \{ 1 \}$, we must show that the
image of $\phi$ in $\Gamma$ is non-trivial.
We do this by showing that, with respect to the given action of $\Gamma$,
$\phi$ does not act as the identity.
The only situation in question is when $\phi$ has normal form $\phi =
\phi_{i_{1}} \cdots \phi_{i_{r}}$ of type 3 from \S \ref{amalgsect}.

Fix any $i \neq i_{r}$.
Then
$$\phi_{i_{r}}(\Theta_{i}) \subset \Theta_{i_{r}}$$
by property 4 of a proper interactive $P$-tuple.
Similarly,
$$\phi_{i_{r-1}}(\phi_{i_{r}}(\Theta_{i})) \subset
\phi_{i_{r-1}}(\Theta_{i_{r}}) \subset \Theta_{i_{r-1}}$$
So, repeatedly applying 4 and inducting, we see see that
$$\phi(\Theta_{i}) \subset \Theta_{i_{1}}$$

Now if $i \neq i_{1}$, then we are done since $\phi$ does not act as the
identity.
If $i = i_{1}$ then $\phi(\Theta_{i}) \subset \Theta_{i}$.
However, by property 5 above, $\phi(\Theta_{i}) \neq \Theta_{i}$ and hence
$\phi$ is again not acting as the identity. $\Box$

%%%%%%%%%%%%%%%%%%%%%%%%%%%%%%%%%%%%%%%%%%%%%%%%%%%%%
%%%%%%%%%%%%%%%%%%%%%%%%%%%%%%%%%%%%%%%%%%%%%%%%%%%%%
\subsection{Kleinian group example} \label{kleinsect}
%%%%%%%%%%%%%%%%%%%%%%%%%%%%%%%%%%%%%%%%%%%%%%%%%%%%%
%%%%%%%%%%%%%%%%%%%%%%%%%%%%%%%%%%%%%%%%%%%%%%%%%%%%%

We now describe the motivating example of the Klein-Maskit combination
theorem in more detail.
We will consider the action of $\PSL_{2}({\mathbb C})$ on
$\widehat{\mathbb C} = {\mathbb C} \cup \{ \infty \}$.
Continue to denote the Fuchsian groups $G_{1}, G_{2} < \PSL_{2}({\mathbb
R})$.
Each stabilizes $\widehat{\mathbb R} = {\mathbb R} \cup \{ \infty \}$, and
after conjugating if necessary, we may assume that the maximal parabolic
subgroup of each is
$$G_{0} = \Stab_{G_{1}}(\infty) = \Stab_{G_{2}}(\infty)
= \left\langle \left( \begin{array}{cc}
1 & 1 \\
0 & 1 \\ \end{array} \right) \right\rangle$$
The elements of $G_{0}$ act on ${\mathbb C}$ by translations parallel to
the ${\mathbb R}$-axis.

Each $G_{i}$ also stabilizes upper and lower half planes, ${\mathbb U}$
and ${\mathbb L}$.
Let $H$ denote the union of the two horoballs
$$H = \{ z \in {\mathbb C} \, | \, \Im(z) > 1 \} \cup \{ z \in {\mathbb C}
\, | \, \Im(z) < -1 \} \subset {\mathbb U } \cup {\mathbb L}$$
It is a well known consequence of the J$\o$rgensen inequality that for
every $\phi \in G_{i} \setminus G_{0}$, $i = 1,2$, we have
\begin{equation} \label{kleinhoreqn}
\phi(H) \subset \Theta = \{ z \in {\mathbb C} \, | \, |\Im(z)| \leq 1 \}
\end{equation}

Now let $h \in \Stab_{\PSL_{2}({\mathbb C})}(\infty) \setminus
\PSL_{2}({\mathbb R})$ be any parabolic.
This has the form
$$h = \left( \begin{array}{cc}
1 & \mu \\
0 & 1 \\ \end{array} \right)$$
for some $\mu \in {\mathbb C} \setminus {\mathbb R}$, and acts on
${\mathbb C}$ by translations {\em transverse} to the ${\mathbb R}$-axis.

It follows that there exists $K > 0$ such that for any $k \geq K$, we have
$$h^{k}(\Theta) \subset H$$
(thus, ``sufficiently complicated'' means that it translates a large
distance transverse to the ${\mathbb R}$-axis).
Set $\Theta_{1} = \Theta$ and $\Theta_{2} = h^{k}(\Theta)$.
This implies
$$\Theta_{1} \subset h^{k}(H) \mbox{ and } \Theta_{2} \subset H$$
We leave it as an exercise using this and (\ref{kleinhoreqn}) to verify
that $\Theta_{1},\Theta_{2}$ is a proper interactive pair for $G_{1},
h^{k} G_{2} h^{-k}$.

%%%%%%%%%%%%%%%%%%%%%%%%%%%%%%%%%%%%%%%%%%%%%%%%%%%%%%%%%%%%%%%%%%%%%%%%
%%%%%%%%%%%%%%%%%%%%%%%%%%%%%%%%%%%%%%%%%%%%%%%%%%%%%%%%%%%%%%%%%%%%%%%%
%%%%%%%%%%%%%%%%%%%%%%%%%%%%%%%%%%%%%%%%%%%%%%%%%%%%%%%%%%%%%%%%%%%%%%%%
\section{Topology and combinatorics of surfaces} \label{topcombsurfsect}
%%%%%%%%%%%%%%%%%%%%%%%%%%%%%%%%%%%%%%%%%%%%%%%%%%%%%%%%%%%%%%%%%%%%%%%%
%%%%%%%%%%%%%%%%%%%%%%%%%%%%%%%%%%%%%%%%%%%%%%%%%%%%%%%%%%%%%%%%%%%%%%%%

We will consider connected orientable surfaces of genus $g$ with $n$
punctures or boundary components and say that this has type $(g,n)$
(whenever we encounter a disconnected surface, we will work with its
components).
We will blur the distinction between a puncture and a boundary component
whenever convenient.
Two exceptions to this are (1) whenever we consider a complex structure on
our surface all boundary components will be replaced by punctures, and (2)
whenever we consider essential arcs on our surface (see \S \ref{sccsect}),
all punctures will be replaced by boundary components.
We sometimes write $R_{g,n}$ to denote a surface of type $(g,n)$ and
$R_{g}$ for a closed surface of genus $g$.

A surface with non-negative Euler characteristic or of type $(0,3)$ is
said to be of {\em excluded type}.
In addition to surfaces of excluded type, those with type $(1,1)$ or
$(0,4)$ are often too small for general arguments and definitions to be
valid, and we refer to these as {\em sporadic}.
We will only concern ourselves with surfaces of non-excluded type.

Given a punctured surface $R$ of type $(g,n)$, we denote the closed
surface with the punctures ``filled back in'' by $\overline{R}$.
We view this as a closed surface of genus $g$ with $n$ marked points.

%%%%%%%%%%%%%%%%%%%%%%%%%%%%%%%%%%%%%%%%%%%%%%%%%
%%%%%%%%%%%%%%%%%%%%%%%%%%%%%%%%%%%%%%%%%%%%%%%%%
\subsection{Simple closed curves} \label{sccsect}
%%%%%%%%%%%%%%%%%%%%%%%%%%%%%%%%%%%%%%%%%%%%%%%%%
%%%%%%%%%%%%%%%%%%%%%%%%%%%%%%%%%%%%%%%%%%%%%%%%%

We denote the set of isotopy classes of essential simple closed curves on
$R$ by ${\cal C}^{0}(R)$.
These are (isotopy classes of) simple closed curves which are
homotopically essential in $R$ and not homotopic to any puncture or
boundary component. A {\em multi-curve} is a finite union of
essential simple closed curves which are pairwise disjoint and pairwise
non-parallel.
The isotopy classes of multi-curves will be denoted ${\cal C}(R)$.
We sometimes view the elements of ${\cal C}(R)$ as finite subsets of
${\cal C}^{0}(R)$.
We will generally confuse isotopy classes with representatives whenever
convenient.

Given $A \in {\cal C}(R)$, we let $R \setminus N(A)$ denote the surface
obtained by removing an open tubular neighborhood of a representative of
$A$ from $R$.
That is $R \setminus N(A)$ denotes {\em $R$ cut open along $A$}.
A multi-curve $A$ is said to be {\em sparse} if every component of $R
\setminus N(A)$ has non-excluded type.

We denote the set of isotopy classes of essential simple closed curves and
essential arcs in $R$ by ${\cal A}^{0}(R)$.
For us, an essential arc is an embedding of the pair $(I,\partial I)
\rightarrow (R,\partial R)$ that cannot be homotoped (rel boundary) into a
boundary component.
Multi-arcs are defined analogously to multi-curves, and
we denote the set of isotopy classes of multi-arcs by ${\cal A}(R)$.

Given $A,B \in {\cal A}(R)$ we denote the {\em geometric intersection
number} of $A$ and $B$ by $i(A,B)$.
This is the minimal number of transverse intersection points over all
representatives of $A$ and $B$.\\

\noindent
{\bf Remark.} In our definitions we do {\em not} require that an isotopy
fix the boundary pointwise.

%%%%%%%%%%%%%%%%%%%%%%%%%%%%%%%%%%%%%%%%%%%%%%%%
%%%%%%%%%%%%%%%%%%%%%%%%%%%%%%%%%%%%%%%%%%%%%%%%
\subsection{Measured laminations} \label{mlsect}
%%%%%%%%%%%%%%%%%%%%%%%%%%%%%%%%%%%%%%%%%%%%%%%%
%%%%%%%%%%%%%%%%%%%%%%%%%%%%%%%%%%%%%%%%%%%%%%%%

In what follows, we fix a complete hyperbolic metric of finite area on $R$
(with geodesic boundary).

%%%%%%%%%%%%%%%%%%%%%%%%%%%%%%%%%%%%%%%%%%%%%%%%%%%
\subsubsection{Geodesic laminations} \label{glsect}
%%%%%%%%%%%%%%%%%%%%%%%%%%%%%%%%%%%%%%%%%%%%%%%%%%%

We denote the space of compact {\em geodesic laminations} on $R$ by ${\cal
GL}_{0}(R)$, and we give this the Thurston topology.
In this topology, a sequence of laminations $\{ {\cal L}_{m}
\}_{m=1}^{\infty}$ converges to a lamination ${\cal L}$ if and only if
every leaf of ${\cal L}$ is a limit of leaves of ${\cal L}_{m}$ (this is
weaker than the Hausdorff topology).

Any lamination with all compact leaves is a finite union of simple closed
geodesics, and any such uniquely determines a multi-curve.
Conversely, any multi-curve has a unique geodesic representative, thus we
may identify ${\cal C}(R)$ (and hence also ${\cal C}^{0}(R)$) as a subset
of ${\cal GL}_{0}(R)$.

%%%%%%%%%%%%%%%%%%%%%%%%%%%%%%%%%%%%%%%%%%%%%%%%%%
\subsubsection{Transverse measures} \label{tmsect}
%%%%%%%%%%%%%%%%%%%%%%%%%%%%%%%%%%%%%%%%%%%%%%%%%%

We denote the space of compactly supported {\em measured laminations} by
$\ML_{0}(R)$.
We view an element of $\ML_{0}(R)$ as a compact geodesic lamination along
with an invariant transverse measure of full support.
For any $\lambda \in \ML_{0}(R)$, let $|\lambda|$ denote the underlying
geodesic lamination and $\lambda$ or $d \lambda$ the transverse measure.
We have an action of ${\mathbb R}_{+}$ on $\ML_{0}(R)$ by scaling the
transverse measure.

${\mathbb R}_{+} \times {\cal C}^{0}(R)$ injects into $\ML_{0}(R)$ by
sending $t \cdot a$ to $t$ times the transverse counting measure on the
geodesic representative of $a$.
Given $\lambda$ and $\lambda'$ in $\ML_{0}(R)$, the total measure of the
``product measure'' $d \lambda \times d \lambda'$ is a natural homogeneous
(with respect to the ${\mathbb R}_{+}$ action) extension of the
intersection number function, which we also denote
$$i: \ML_{0}(R) \times \ML_{0}(R) \rightarrow {\mathbb R}$$
We endow $\ML_{0}(R)$ with the smallest topology for which $i$ is
continuous.
With this topology, ${\mathbb R}_{+} \times {\cal C}^{0}(R)$ is dense.

Given a pair $\lambda_{1},\lambda_{2} \in \ML_{0}(R)$, with
$|\lambda_{1}|$ transverse to $|\lambda_{2}|$, we say that $\lambda_{1},
\lambda_{2}$ {\em bind} $R$ if for every $a \in {\cal C}^{0}(R)$, we have
$i(a,\lambda_{1}) + i(a,\lambda_{2}) > 0$.
Equivalently, all complementary regions of $|\lambda_{1}| \cup
|\lambda_{2}|$ are disks, punctured disks, or half-open annuli.

The quotient by ${\mathbb R}_{+}$
$$\prj : \ML_{0}(R) \rightarrow \prj\ML_{0}(R)$$
is the space of {\em projective measured laminations} on $R$.
This is given the quotient topology and we denote $\prj(\lambda)$ by
$[\lambda]$, for any $\lambda \in \ML_{0}(R)$.
The forgetful map to ${\cal GL}_{0}(R)$ obtained by forgetting the
transverse measure factors through the projectivization.

An important point for us is the following result of Thurston \cite{Tn}
Proposition 8.10.3.

\begin{proposition}
The forgetful map
$$\prj\ML_{0}(R) \rightarrow {\cal GL}_{0}(R)$$
is continuous.
\end{proposition}

\noindent
{\bf Remark.} Although it appears that the spaces we have constructed
depend on the choice of hyperbolic metric, any two hyperbolic metrics give
rise to spaces which can be canonically identified.

%%%%%%%%%%%%%%%%%%%%%%%%%%%%%%%%%%%%%%%%%%%%%%%%%%%%%%%%%%%%%%%%%%
\subsubsection{Zero locus of a measured lamination} \label{zlsect}
%%%%%%%%%%%%%%%%%%%%%%%%%%%%%%%%%%%%%%%%%%%%%%%%%%%%%%%%%%%%%%%%%%

For any $\lambda \in \ML_{0}(R)$, $|\lambda|$ is a finite union of its
components $|\lambda| = |\lambda_{1}| \cup \cdots \cup |\lambda_{N}|$.
This decomposes $\lambda$ as a sum of sub-measured laminations $\lambda =
\sum_{j=1}^{N} \lambda_{j}$.
With this notation, we define
$$Z(\lambda) = \{ \mu \in \ML_{0}(R) \, | \, i(\mu,\lambda) = 0 \} \quad
\mbox{ and } \quad Z'(\lambda) = \bigcup_{j=1}^{N} Z(\lambda_{j})$$

These sets depend only on $|\lambda|$.
In particular, if $A$ is a multi-curve, we may define $Z(A)$ and $Z'(A)$
by arbitrarily prescribing a transverse measure of full support on the
geodesic representative of $A$.

By an abuse of notation, we refer to the images of $Z(\lambda)$ and
$Z'(\lambda)$ in $\prj\ML_{0}(R)$ and ${\cal GL}_{0}(R)$ by the same
names.
Note that if $\mu \in Z(\lambda)$ and $\mu' \in \ML_{0}(R)$ has $|\mu| =
|\mu'|$, then $\mu' \in Z(\lambda)$, and similarly for $Z'(\lambda)$.
In particular, we may pass back and forth between the different sets
representing $Z(\lambda)$ and $Z'(\lambda)$ by taking images and preimages
without introducing any new points.

As a further abuse of notation, we will denote the intersection of
$Z(A),Z'(A) \subset {\cal GL}_{0}(R)$ with ${\cal C}(R)$ by the same
names.
It will be clear from the context where these subsets are residing.

%%%%%%%%%%%%%%%%%%%%%%%%%%%%%%%%%%%%%%%%%%%%%%%%%%%%%%%
\subsubsection{Singular foliations} \label{singfolsect}
%%%%%%%%%%%%%%%%%%%%%%%%%%%%%%%%%%%%%%%%%%%%%%%%%%%%%%%

The elements of $\ML_{0}(R)$ can also be realized as singular foliations
on $S$ with invariant transverse measures of full support.
The underlying singular foliation is unique up to isotopy and Whitehead
equivalence.
The corresponding geodesic lamination is obtained by straightening all
leaves to geodesics.
We will also write $|\lambda|$ to denote the underlying (isotopy/Whitehead
equivalence class of) singular foliation.
When it is clear from the context, we will not explicitly mention which of
these we are referring to.

A foliation is {\em finite} if it has only compact leaves.
Given $\lambda \in \ML_{0}(R)$, the foliation $|\lambda|$ is finite if and
only if the geodesic lamination $|\lambda|$ represents a multi-curve.
A finite foliation decomposes $R$ into a finite union of compact singular
leaves and open annuli foliated by its nonsingular leaves.
We refer to this as an {\em annular decomposition}.
The union of the cores of the annuli will be called the {\em core} of the
annular decomposition, and this is precisely the associated multi-curve
for the finite foliation.
An annular decomposition is {\em sparse} if its core is.

%%%%%%%%%%%%%%%%%%%%%%%%%%%%%%%%%%%%%%%%%%%%%
%%%%%%%%%%%%%%%%%%%%%%%%%%%%%%%%%%%%%%%%%%%%%
\subsection{The arc complex} \label{ccsect}
%%%%%%%%%%%%%%%%%%%%%%%%%%%%%%%%%%%%%%%%%%%%%
%%%%%%%%%%%%%%%%%%%%%%%%%%%%%%%%%%%%%%%%%%%%%

The set ${\cal A}^{0}(R)$ can be naturally identified with the vertex set
of a simplicial complex ${\cal A}(R)$, called the {\em arc complex} of
$R$.
A set of distinct elements $a_{0},...,a_{p} \in {\cal A}^{0}$ spans a
$p$-simplex if and only if the curves and arcs can be realized on the
surface pairwise disjointly.
Thus, the points of the complex ${\cal A}(R)$ (that is, its simplices) are
precisely the multi-arcs.
We are therefore justified in using ${\cal A}(R)$ to refer to either of
these objects.

We denote the corresponding polyhedron $|{\cal A}|(R)$.
This is not given the classical topology for the polyhedron of a
simplicial complex, but rather the metric topology obtained by declaring
each simplex to be a regular Euclidean simplex with all edge lengths equal
to one.
In particular, we obtain a metric on ${\cal A}^{0}(R)$ by viewing it as
the set of vertices of the metric graph $|{\cal A}^{1}|(R)$.
Whenever we refer to the distance between elements of ${\cal A}^{0}(R)$,
it will be understood that it is meant with respect to this metrics.

%%%%%%%%%%%%%%%%%%%%%%%%%%%%%%%%%%%%%%%%%%%%%%%%%%%%%%
%%%%%%%%%%%%%%%%%%%%%%%%%%%%%%%%%%%%%%%%%%%%%%%%%%%%%%
\subsection{Subsurface projections} \label{subprosect}
%%%%%%%%%%%%%%%%%%%%%%%%%%%%%%%%%%%%%%%%%%%%%%%%%%%%%%
%%%%%%%%%%%%%%%%%%%%%%%%%%%%%%%%%%%%%%%%%%%%%%%%%%%%%%

A useful tool is the subsurface projection.
If $S \subset R$ is a connected incompressible subsurface of non-excluded
type, then we have the {\em subsurface projection}
$$\pi_{S}:\prj\ML_{0}(R) \setminus Z(\partial S) \rightarrow {\cal
A}^{0}(S)$$
obtained as follows.
We fix a hyperbolic metric on $R$ and homotope the embedding of $S$ to
have geodesic boundary (this may make the boundary no longer embedded, but
the interior remains embedded).
For $[\lambda] \in \prj\ML_{0}(R) \setminus Z(\partial S)$, we consider
the geodesic lamination $|\lambda|$ and take its intersection with $S$
(less any components of $\partial S$).
The result is a disjoint union of essential arcs and a geodesic
laminations in $S$ (we are using the fact that the transverse measure has
full support to guarantee that $|\lambda|$ has no isolated leaves which
are not simple closed curves).
Since $[\lambda]$ is not in $Z(\partial S)$, some leaf of $|\lambda|$
intersects $\partial S$ non-trivially, and we take $\pi_{S}([\lambda])$ to
be the isotopy class of one of these arcs.
The collection of essential arcs and the geodesic lamination in $S$
obtained in this way depend only on the underlying geodesic lamination
$|\lambda|$, and not the projective class of transverse measure.
There is no canonical choice for this arc, however any two such arcs are
disjoint and hence have distance one from each other in ${\cal A}^{0}(S)$.

The use of the hyperbolic structure in this construction is merely a
convenient way of arranging that all intersection points of leaves of
$|\lambda|$ with $\partial S$ are essential, and so the construction may
be carried out assuming only this.

The subsurface projection is ``coarsely continuous''.

\begin{proposition} \label{coarsecontprop}
If $\{ [\lambda_{m}] \}_{m=1}^{\infty} \subset \prj\ML_{0}(R) \setminus
Z(\partial S)$ converges to $[\lambda] \in \prj\ML_{0}(R) \setminus
Z(\partial S)$ then there exists $N > 0$, such that for all $m \geq N$
$$d(\pi_{S}([\lambda]),\pi_{S}([\lambda_{m}])) \leq 1$$
\end{proposition}

\noindent
{\em Proof.} Assume the setup from the statement of the proposition, and
let $l$ be any leaf of $|\lambda|$ which nontrivially intersects
$Z(\partial S)$.
If we let $l^{0}$ be an arc of $l \cap S$, then this is an allowable
choice for $l^{0} = \pi_{S}([\lambda])$.
Since the forgetful map to ${\cal GL}_{0}(R)$ is continuous, there are
leaves $l_{m}$ of $|\lambda_{m}|$ which converge to $l$.
Thus there are arcs $l^{0}_{m}$ of $l_{m} \cap S$ which converge to
$l^{0}$.
In particular, there exists some $N> 0 $ such that for all $m \geq N$, the
arcs $l^{0}_{m}$ are isotopic in $S$ to $l^{0}$.

Clearly the same thing holds for any choice $l^{0}$ for
$\pi_{S}([\lambda])$.
So, since any choice of arc for $\Pi_{S}([\lambda_{m}])$ will be disjoint
from the $l^{0}_{m}$ considered above, we obtain the required bound.
$\Box$ \\

\noindent
{\bf Remark.} As is noted there, Proposition \ref{coarsecontprop} is the
intuition behind the Bounded Geodesic Image Theorem of \cite{MM2}.\\

Since $\pi_{S}(\lambda)$ depends only on $|\lambda|$, it makes sense to
define $\pi_{S}(A)$ for $A \in {\cal C}(R) \setminus Z(\partial S)$.
However, we will adhere to the following

\begin{convention} \label{importantconv}
When considering multi-curves, the domain for the subsurface projections
will be taken to be
$${\cal C}(R) \setminus Z'(\partial S)$$
instead of ${\cal C}(R) \setminus Z(\partial S)$
\end{convention}

%%%%%%%%%%%%%%%%%%%%%%%%%%%%%%%%%%%%%%%%%%%%%%%%
%%%%%%%%%%%%%%%%%%%%%%%%%%%%%%%%%%%%%%%%%%%%%%%%
%%%%%%%%%%%%%%%%%%%%%%%%%%%%%%%%%%%%%%%%%%%%%%%%
\section{Quadratic differentials} \label{qdsect}
%%%%%%%%%%%%%%%%%%%%%%%%%%%%%%%%%%%%%%%%%%%%%%%%
%%%%%%%%%%%%%%%%%%%%%%%%%%%%%%%%%%%%%%%%%%%%%%%%
%%%%%%%%%%%%%%%%%%%%%%%%%%%%%%%%%%%%%%%%%%%%%%%%

The {\em Teichm\"uller space} of $R$, denoted ${\cal T}(R)$, is the space
of isotopy classes complex structures of finite type on $R$.
For any $X \in {\cal T}(R)$ let $A_{2}(X)$ denote the vector space of
integrable holomorphic quadratic differentials on $R$ with respect to $X$.
Varying $X$ over ${\cal T}(R)$ we obtain a vector bundle over ${\cal
T}(R)$ with fiber over $X$ being $A_{2}(X)$.
We denote the complement of the zero section of this bundle as $A_{2}{\cal
T}(R)$.
We refer to a point of $A_{2}{\cal T}(R)$ simply by the quadratic
differential $q$, with the underlying complex structure $X \in {\cal
T}(R)$ implicit.

%%%%%%%%%%%%%%%%%%%%%%%%%%%%%%%%%%%%%%%%%%%%%%%%
%%%%%%%%%%%%%%%%%%%%%%%%%%%%%%%%%%%%%%%%%%%%%%%%
\subsection{Metric structure} \label{metricsect}
%%%%%%%%%%%%%%%%%%%%%%%%%%%%%%%%%%%%%%%%%%%%%%%%
%%%%%%%%%%%%%%%%%%%%%%%%%%%%%%%%%%%%%%%%%%%%%%%%

Any point $q \in A_{2}{\cal T}(R)$ determines, and is determined by, a
Euclidean cone metric on $\overline{R}$ with some additional structure
which we now describe.
There is an atlas of local coordinate charts $\{ \zeta_{j} : U_{j}
\rightarrow {\mathbb C} \}$ on the complement of the cone points with the
following properties:
\begin{enumerate}
\item for each $j,k$ and on each component of $U_{j} \cap U_{k}$, there is
a $w \in {\mathbb C}$ such that
$$\zeta_{j} = \pm \zeta_{k} + w$$
\item the metric on $U_{j}$ is the pullback of the Euclidean metric on
${\mathbb C}$ by $\zeta_{j}$.
\end{enumerate}
It follows that cone points have cone angle $m \pi$ for some $m \in
{\mathbb Z}_{\geq 1}$.
We further require that at any cone point which is not a marked point, $m
\geq 3$.

Any local coordinate as above is called a $q$-coordinate.
We refer to the set of cone points and marked points as the singularities
of $q$.
Because of condition (1) on the transition functions, this structure on
$\overline{R}$ is sometimes called a {\em half-translation structure}.

Any embedded geodesic segment in the $q$ metric connecting a pair of
singular points and having no singular point in its interior is called a
{\em saddle connection}.
A simple closed geodesic which misses the singular locus will be called a
{\em nonsingular closed geodesic}.
There is a dichotomy for simple closed curves $a \in {\cal C}^{0}(R)$ with
respect to $q$: either there is an annulus in $R$ foliated by nonsingular
closed geodesics all in the isotopy class of $a$, or else there is a
unique closed geodesic in $\overline{R}$ homotopic to $a$ made up of
saddle connections.

The geodesic representative of $a \in {\cal C}^{0}(R)$ when $R$ has
punctures requires a few additional comments.
If we take a sequence of representatives of $a$ in $R$ for which the
lengths converge to the infimum of all lengths of representatives, then a
subsequence of these representatives will converge to a geodesic {\em in
$\overline{R}$}, and this is the geodesic representative of $a$.
In particular, arbitrarily close to the geodesic representative we can
find simple closed curves in $R$ representing $a$.

%%%%%%%%%%%%%%%%%%%%%%%%%%%%%%%%%%%%%%%%%
%%%%%%%%%%%%%%%%%%%%%%%%%%%%%%%%%%%%%%%%%
\subsection{$\ML_{0}(q)$} \label{mlqsect}
%%%%%%%%%%%%%%%%%%%%%%%%%%%%%%%%%%%%%%%%%
%%%%%%%%%%%%%%%%%%%%%%%%%%%%%%%%%%%%%%%%%

For any constant, nonzero $1$-form $a dx + b dy$ on ${\mathbb C}$ (so $a,b
\in {\mathbb R}$), the measured foliation $|a dx + b dy|$ is invariant
under the transition functions for the atlas of $q$-coordinates.
Therefore, it pulls back to a measured foliation on the complement of the
cone points.
This extends to a singular measured foliation on $R$ and $\overline{R}$.

The space of all such measured foliations is denoted $\ML_{0}(q)$ and is
homeomorphic to $({\mathbb R}^{2} \setminus \{ 0 \}) / \pm I$.
We view $\ML_{0}(q)$ as a subspace of $\ML_{0}(R)$ and we write
$\prj\ML_{0}(q) \subset \prj\ML_{0}(R)$ to denote the projectivization.

We note that the leaves of every element of $\ML_{0}(q)$ are geodesic in
the Euclidean cone metric, with a finite number of {\em critical leaves}
containing one or two singular points of the metric as endpoints.
When there are two singular points, this is just a saddle connection.

When $\lambda \in \ML_{0}(q)$ is a finite foliation, it defines an annular
decomposition of $R$ (and so also a core for this annular decomposition---see \S \ref{singfolsect}).
We call this is an {\em annular decomposition} for $q$.

\begin{lemma} \label{allinterlem}
If $A_{0},A_{1}$ are cores of two distinct annular decompositions for $q$,
then $A_{0} \not \in Z'(A_{1})$.
\end{lemma}

\noindent
{\em Proof.} For any two distinct $[\lambda_{0}],[\lambda_{1}] \in
\prj\ML_{0}(q)$, it is easy to see that $\lambda_{0},\lambda_{1}$ bind
$R$. In particular $A_{0},A_{1}$ bind $R$.
Every component $a \subset A_{1}$ must therefore intersect $A_{0}$
(otherwise $a$ would have $i(a,A_{0}) + i(a,A_{1}) = 0$). $\Box$

%%%%%%%%%%%%%%%%%%%%%%%%%%%%%%%%%%%%%%%%%%%%%
%%%%%%%%%%%%%%%%%%%%%%%%%%%%%%%%%%%%%%%%%%%%%
%%%%%%%%%%%%%%%%%%%%%%%%%%%%%%%%%%%%%%%%%%%%%
\section{Mapping class group} \label{mcgsect}
%%%%%%%%%%%%%%%%%%%%%%%%%%%%%%%%%%%%%%%%%%%%%
%%%%%%%%%%%%%%%%%%%%%%%%%%%%%%%%%%%%%%%%%%%%%
%%%%%%%%%%%%%%%%%%%%%%%%%%%%%%%%%%%%%%%%%%%%%

%%%%%%%%%%%%%%%%%%%%%%%%%%%%%%%%%%%%%%%
%%%%%%%%%%%%%%%%%%%%%%%%%%%%%%%%%%%%%%%
\subsection{Actions} \label{mcgactsect}
%%%%%%%%%%%%%%%%%%%%%%%%%%%%%%%%%%%%%%%
%%%%%%%%%%%%%%%%%%%%%%%%%%%%%%%%%%%%%%%

$\Mod(R)$ acts on all the spaces associated to $R$ we have discussed so
far: ${\cal A}^{0}(R)$, ${\cal C}(R)$, $\ML_{0}(R)$, $\prj \ML_{0}(R)$, and
${\cal GL}_{0}(R)$.
These actions respect the various natural maps between the spaces (e.g.
the inclusion ${\cal C}^{0}(R) \subset \ML_{0}(R)$, and the projection
$\ML(R) \rightarrow \prj \ML(R)$).
These actions also preserves the intersection number $i$ (on the spaces
where it is defined) as well as the metric on ${\cal A}^{0}(R)$.

%%%%%%%%%%%%%%%%%%%%%%%%%%%%%%%%%%%%%%%%%%%%%%%%%%%%%%%%%%%%%%
%%%%%%%%%%%%%%%%%%%%%%%%%%%%%%%%%%%%%%%%%%%%%%%%%%%%%%%%%%%%%%
\subsection{Nielsen-Thurston classification} \label{classsect}
%%%%%%%%%%%%%%%%%%%%%%%%%%%%%%%%%%%%%%%%%%%%%%%%%%%%%%%%%%%%%%
%%%%%%%%%%%%%%%%%%%%%%%%%%%%%%%%%%%%%%%%%%%%%%%%%%%%%%%%%%%%%%

An automorphism $\phi \in \Mod(R) \setminus \{ 1 \}$ is said to be {\em
reducible} if there exists a multi-curve $A \in {\cal C}(R)$ invariant by
$\phi$.
In this case, write $S = R \setminus N(A)$ and denote the components of
$S$ by $S_{1},...,S_{M}$.
We thus have an induced automorphism
$$\tilde{\phi}:S \rightarrow S$$
Some power $\tilde{\phi}^{r}$, preserves the components of this surface
and hence gives rise to automorphisms
$$\tilde{\phi}^{r}_{j}:S_{j} \rightarrow S_{j}$$
The multi-curve $A$ is called a {\em reducing system} for $\phi$ and the
maps $\tilde{\phi}^{r}_{j}$ are called the {\em components} of
$\tilde{\phi}^{r}$ (or of $\phi^{r}$).
The process just described for a reducible automorphism is called a {\em
reduction} of $\phi$ along $A$.

An example of a reducible automorphism is a Dehn twist $T_{a}$ along some
simple closed curve $a \in {\cal C}^{0}(R)$.
More generally, given a multi-curve $A$, a composition of Dehn twists
about the components of $A$ defines a reducible automorphism called a {\em
multi-twist}.
If all the Dehn twists are positive, it is called a {\em positive
multi-twist}.\\

An element $\phi \in \Mod(R)$ of infinite order is called {\em
pseudo-Anosov} if there exists a pair $\lambda_{s},\lambda_{u} \in
\ML_{0}(R)$ called the {\em stable} and {\em unstable} laminations,
respectively, which bind $R$ and such that $\{ [\lambda_{s}],
[\lambda_{u}] \}$ is invariant by $\phi$.  

\begin{theorem} {\bf (Nielsen-Thurston)} \label{autclassthrm}
Given any element $\phi \in \Mod(R)$, one of the following holds
\begin{enumerate}
\item \label{foitem} $\phi$ has finite order,
\item \label{reditem}$\phi$ is reducible, or
\item \label{paitem} $\phi$ is pseudo-Anosov.
\end{enumerate}
\end{theorem}

If $\phi$ is reducible, there is a unique reducing system $A$ with the
property that, after passing to an appropriate power $\phi^{r}$, each of
the components $\tilde{\phi}^{r}_{j}$ are either pseudo-Anosov or the
identity, and any sub-multi-curve of $A$ fails to have this property.
We refer to this multi-curve as the {\em canonical reducing system} for
$\phi$ and those components of (a power of) $\phi$ which are pseudo-Anosov
are called the {\em pseudo-Anosov components}.

Whenever we perform a reduction, it will always be assumed to be done
along the canonical reducing system.
We say that $\phi$ is {\em pure} if it does not permute the components of
$S = R \setminus N(A)$ or the components of $A$ and each component of
$\phi$ is either pseudo-Anosov or the identity.
Any reducible automorphism has a power which is pure.\\

%%%%%%%%%%%%%%%%%%%%%%%%%%%%%%%%%%%%%%%%%%%%%%%%%%%%%%%%%%%%%%%%%%%%%%%%%%%
%%%%%%%%%%%%%%%%%%%%%%%%%%%%%%%%%%%%%%%%%%%%%%%%%%%%%%%%%%%%%%%%%%%%%%%%%%%
\subsection{Identifying the type of an automorphism} \label{identpasect}
%%%%%%%%%%%%%%%%%%%%%%%%%%%%%%%%%%%%%%%%%%%%%%%%%%%%%%%%%%%%%%%%%%%%%%%%%%%
%%%%%%%%%%%%%%%%%%%%%%%%%%%%%%%%%%%%%%%%%%%%%%%%%%%%%%%%%%%%%%%%%%%%%%%%%%%

Given $\phi \in \Mod(R)$, one can decide which of the three types from
Theorem \ref{autclassthrm} $\phi$ falls into by considering the action on
$\prj\ML_{0}(R)$.

In case \ref{foitem} of Theorem \ref{autclassthrm}, there is nothing to
say (some power of $\phi$ acts as the identity).
In case \ref{paitem}, the two points $[\lambda_{s}]$ and $[\lambda_{u}]$
are attracting and repelling fixed points for the action, respectively.
In particular for every $[\lambda] \in \prj\ML_{0}(R) \setminus \{
[\lambda_{s}], [\lambda_{u}] \}$, we have
$$\phi^{n}([\lambda]) \rightarrow [\lambda_{s}] \mbox{ and }
\phi^{-n}([\lambda]) \rightarrow [\lambda_{u}]$$
as $n \rightarrow \infty$.

For simplicity, we only describe the reducible case when $\phi$ is pure.
Let $A$ denote a canonical reducing system for $\phi$.
Let $\lambda_{s,1},...,\lambda_{s,M},\lambda_{u,1},...,\lambda_{u,M}$
denote the stable and unstable measured laminations for the pseudo-Anosov
components of $\phi$, which we view as elements of $\ML_{0}(R)$.
>From this data we obtain two laminations
$$\psi_{s} = A \cup |\lambda_{s,1}| \cup \cdots \cup |\lambda_{s,M}|
\mbox{ and } \psi_{u} = A \cup |\lambda_{u,1}| \cup \cdots \cup
|\lambda_{u,M}|$$
These define subsets $\Psi_{s}$ and $\Psi_{u}$ of $\ML_{0}(R)$ consisting
of all measured laminations with measures supported on (non-empty)
sublaminations of $\psi_{s}$ and $\psi_{u}$, respectively.
We will need the following result of Ivanov and McCarthy, see \cite{Iv} Theorem
A1.

\begin{theorem} \label{autotypethrm}
With the notation as above, for every $\lambda \in \ML_{0}(R) \setminus
Z(\Psi_{s}) \cup Z(\Psi_{u})$
$$\lim_{m \rightarrow \infty} \phi^{m}([\lambda]) \in \prj\Psi_{s} \mbox{
and } \lim_{m \rightarrow -\infty} \phi^{m}([\lambda]) \in \prj\Psi_{u}$$
\end{theorem}

Note that if $\lambda \in Z(\Psi_{s}) \cap Z(\Psi_{u})$, then $\phi$ fixes
$\lambda$.

%%%%%%%%%%%%%%%%%%%%%%%%%%%%%%%%%%%%%%%%%%%%%%%%%%%%%%%%%%%%%%%%%%%%%%%%%%
%%%%%%%%%%%%%%%%%%%%%%%%%%%%%%%%%%%%%%%%%%%%%%%%%%%%%%%%%%%%%%%%%%%%%%%%%%
\subsection{Pseudo-Anosov automorphisms and ${\cal A}^{0}(R)$}
\label{passect}
%%%%%%%%%%%%%%%%%%%%%%%%%%%%%%%%%%%%%%%%%%%%%%%%%%%%%%%%%%%%%%%%%%%%%%%%%%
%%%%%%%%%%%%%%%%%%%%%%%%%%%%%%%%%%%%%%%%%%%%%%%%%%%%%%%%%%%%%%%%%%%%%%%%%%

The following theorem of Masur and Minsky is a strengthening of the fact
that a pseudo-Anosov automorphism cannot fix any essential simple closed
curve or arc.

\begin{theorem} \label{mmpathrm} {\bf (Masur-Minsky)}
There exists $c > 0$ (depending only on $R$), such that for any
pseudo-Anosov automorphism $\phi \in \Mod(R)$ the action of $\phi$ on
${\cal A}^{0}(R)$ satisfies
$$d(\phi^{k}(a),a) \geq c |k|$$
for every $a \in {\cal A}^{0}(R)$.
\end{theorem}

\noindent
{\bf Remark.} Although this is stated for the {\em curve complex} instead
of the arc complex, the two spaces are quasi-isometric by a
$\Mod(R)$-equivariant quasi-isometry (in the sporadic cases, the curve
complex is replaced by the Farey graph).
This easily implies the version we need.

%%%%%%%%%%%%%%%%%%%%%%%%%%%%%%%%%%%%%%%%%%%
%%%%%%%%%%%%%%%%%%%%%%%%%%%%%%%%%%%%%%%%%%%
\subsection{Veech groups} \label{veechsect}
%%%%%%%%%%%%%%%%%%%%%%%%%%%%%%%%%%%%%%%%%%%
%%%%%%%%%%%%%%%%%%%%%%%%%%%%%%%%%%%%%%%%%%%

Let $\Aff^{+}(q)$ denote the subgroup of $\Mod(R)$ consisting of all
automorphisms with representatives which are affine with respect to the
$q$-metric.
We will refer to any subgroup $G(q) < \Aff^{+}(q)$ as a {\em Veech group}
for $q$. Note this is non-standard notation, but is convenient for
our purposes.

Taking derivatives in $q$-coordinates (with respect to the standard basis
of ${\mathbb C} \cong {\mathbb R}^{2}$) we obtain a homomorphism
$$D:\Aff^{+}(q) \rightarrow \PSL_{2}({\mathbb R})$$
This gives us the short exact sequence shown below, and 
in which $\PSL(q)$ is the image group of the homomorphism $D$ above,
and $\Aut(q)$ is the subgroup of $\Aff^{+}(q)$ preserving $q$.

\begin{equation} \label{seseqn}
1 \rightarrow \Aut(q) \rightarrow \Aff^{+}(q) \rightarrow \PSL(q)
\rightarrow 1
\end{equation}

The next theorem is well known.
It can be derived from Bers proof of Theorem \ref{autclassthrm}, as is
done by Kra in \cite{Kra}.
A proof is also given by Veech in \cite{V1} and Thurston in \cite{T}.

\begin{theorem} \label{affactthrm}
$D$ is a discrete representation with finite kernel.
For $\phi \in \Aff^{+}(q) \setminus \{ 1 \}$:
\begin{enumerate}
\item \label{elliptic} $\phi$ has finite order if and only if $D\phi$ is
elliptic or $D\phi=1$,
\item \label{parabolic} $\phi$ has a power that is a positive multi-twist
if and only if $D\phi$ is parabolic, and
\item \label{hyperbolic} $\phi$ is pseudo-Anosov if and only if $D\phi$ is
hyperbolic.
\end{enumerate}
\end{theorem}

The proof boils down to the following elementary observation.
Given $\phi \in \Aff^{+}(q)$, either $D \phi$ preserves a quadratic form
(case \ref{elliptic}), has a single eigenspace (case \ref{parabolic}), or
a pair of distinct eigenspaces (case \ref{hyperbolic}).
Then, either $\phi$ preserves a metric affine equivalent to the $q$-metric
(case \ref{elliptic}) or a (pair of) foliation(s) everywhere tangent to
the eigenspace(s) (cases \ref{parabolic} and \ref{hyperbolic}).

In case \ref{parabolic}, the invariant foliation is finite and so defines
an annular decomposition with core $A_{0}$.
Some power of $\phi$ is a positive multi-twist about $A_{0}$.
The stabilizer of $A_{0}$ in $\Aff^{+}(q)$ contains $\langle \phi \rangle$
with finite index.
For any Veech group $G(q) < \Aff^{+}(q)$ we denote the stabilizer of
$A_{0}$ in $G(q)$ by $G_{0}(q)$.\\

We say that a subgroup of $\Aff^{+}(R)$ has {\em elliptic} or {\em
parabolic type} if its image under $D$ is an elementary subgroup with
elliptic or parabolic type, respectively.
In particular, $G_{0}(q)$ is of parabolic type.

%%%%%%%%%%%%%%%%%%%%%%%%%%%%%%%%%%%%
%%%%%%%%%%%%%%%%%%%%%%%%%%%%%%%%%%%%
\subsection{Examples} \label{exsect}
%%%%%%%%%%%%%%%%%%%%%%%%%%%%%%%%%%%%
%%%%%%%%%%%%%%%%%%%%%%%%%%%%%%%%%%%%

\begin{figure}[htb]
\centerline{}
\begin{center}
\ \psfig{file=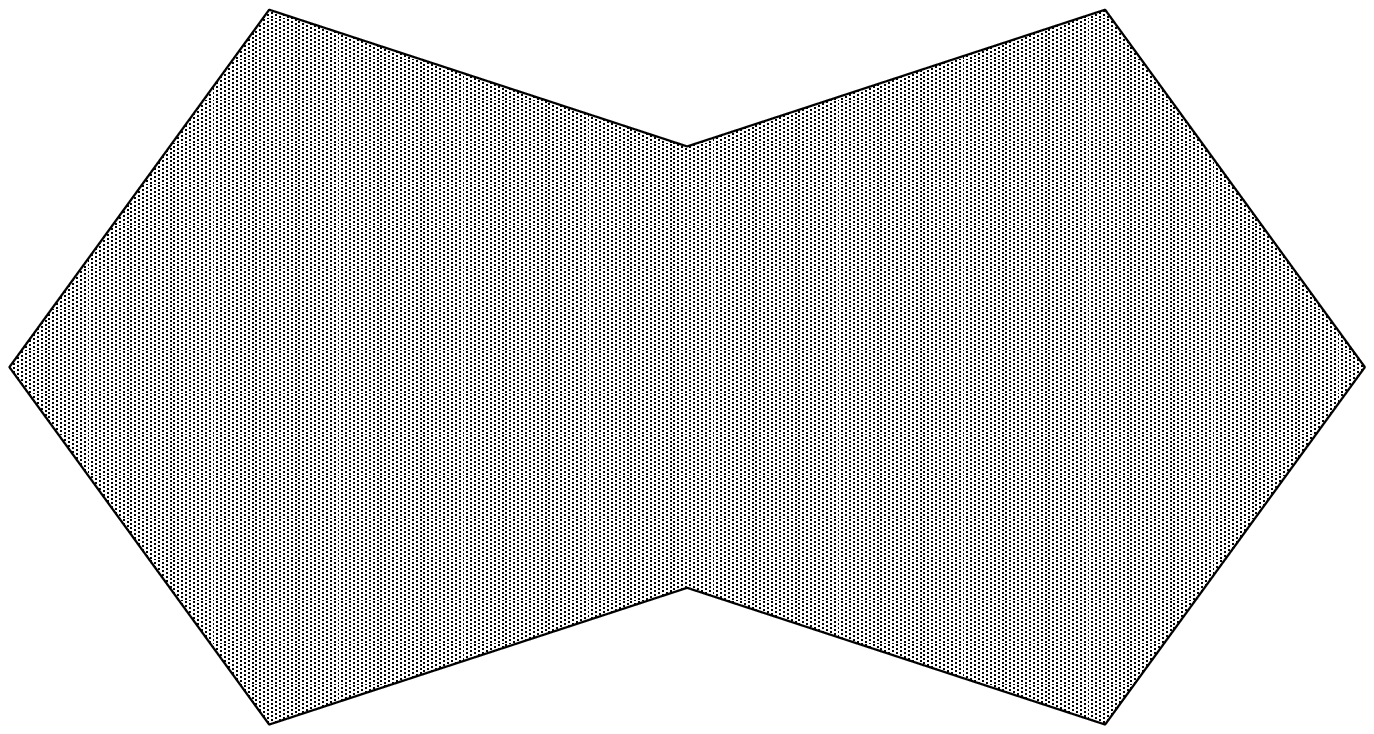,height=1truein}
\caption{$\Delta_{2}$} \label{pentagonsfig}
\end{center}
\end{figure}

The following construction is due to Veech.
Given $g \geq 2$, let $\Delta_{g}$ be the non-convex polygon obtained as
the union of two regular $2g+1$-gons in the Euclidean plane which meet
along an edge and have disjoint interiors (see Figure \ref{pentagonsfig}
for the case $g = 2$).

Let $R_{g}$ denote the surface obtained by gluing opposite sides of
$\Delta_{g}$ by translations.
This is easily seen to be a closed surface of genus $g$.
The Euclidean metric on the interior of $\Delta_{g}$ is the restriction of
a Euclidean cone metric on $R_{g}$, and we can find local coordinates
defining a quadratic differential on $R_{g}$ compatible with this metric
(as in \ref{metricsect}).
We denote this quadratic differential $\xi_{g} \in A_{2}{\cal T}(R_{g})$.

\begin{theorem} {\bf (Veech)} \label{veechexthrm}
$\PSL(\xi_{g})$ is a lattice isomorphic to a triangle group of type
$(2,2g+1,\infty)$.
The single primitive parabolic conjugacy class is the image by $D$ of the
conjugacy class of a positive multi-twist in $\Aff^{+}(\xi_{g})$ by $D$
about a sparse multi-curve $A_{0} \in {\cal C}(R_{g})$.
\end{theorem}

See \cite{V1} for a proof.\\

$\PSL(\xi_{g})$ is generated by a pair of elliptic elements $\gamma_{0}$
and $\gamma_{1}$ with orders $2$ and $2g+1$, respectively.
Moreover, these can be chosen so that $\gamma_{0} \gamma_{1}$ generates
the parabolic subgroup.

There is an epimorphism
$$\nu:\PSL(\xi_{g}) \rightarrow {\mathbb Z}/2{\mathbb Z} \oplus {\mathbb
Z}/(2g+1){\mathbb Z}$$
with $\nu(\gamma_{0})$ and $\nu(\gamma_{1})$ generating the first and
second factors, respectively.
It follows that $\ker(\nu)$ is a free group.
Moreover, since $\nu(\gamma_{0} \gamma_{1})$ generates the range, we see
that $\ker(\nu)$ has exactly one conjugacy class of parabolics.
In fact, $\ker(\nu) \cong \pi_{1}(S_{g,1})$ with the parabolic
$(\gamma_{0}\gamma_{1})^{2(2g+1)}$ representing a peripheral element,
where $S_{g,1}$ is a surface of type $(g,1)$.

The sequence (\ref{seseqn}) restricts to a short exact sequence
$$1 \rightarrow D^{-1}(\ker(\nu)) \cap \Aut(\xi_{g}) \rightarrow
D^{-1}(\ker(\nu)) \rightarrow \ker(\nu) \rightarrow 1$$
which splits since $\ker(\nu)$ is free.
Denoting the image of $\ker(\nu)$ under the splitting isomorphism by
$G(\xi_{g}) < \Aff^{+}(\xi_{g})$ we have

\begin{corollary} \label{onepunctcor}
$G(\xi_{g}) < \Mod(R_{g})$ is isomorphic to $\pi_{1}(S_{g,1})$ with
every non-peripheral conjugacy class represented by a pseudo-Anosov
automorphism and the peripheral subgroup generated by a positive
multi-twist about the sparse multi-curve $A_{0} \in {\cal C}(R_{g})$.
\end{corollary}

\noindent
{\em Proof.} This follows from Theorem \ref{veechexthrm} and the comments
above by appealing to Theorem \ref{affactthrm}. $\Box$

%%%%%%%%%%%%%%%%%%%%%%%%%%%%%%%%%%%%%%%%%%%%%%%%%%%
%%%%%%%%%%%%%%%%%%%%%%%%%%%%%%%%%%%%%%%%%%%%%%%%%%%
%%%%%%%%%%%%%%%%%%%%%%%%%%%%%%%%%%%%%%%%%%%%%%%%%%%
\section{The combination theorem} \label{combosect}
%%%%%%%%%%%%%%%%%%%%%%%%%%%%%%%%%%%%%%%%%%%%%%%%%%%
%%%%%%%%%%%%%%%%%%%%%%%%%%%%%%%%%%%%%%%%%%%%%%%%%%%
%%%%%%%%%%%%%%%%%%%%%%%%%%%%%%%%%%%%%%%%%%%%%%%%%%%

{\em For the remainder of this paper, we let $R$ denote a non-sporadic
surface of type $(g,n)$.}\\

We first define the groups to which our theorem applies.
Let $q_{1},....,q_{P} \in A_{2}{\cal T}(R)$ and suppose that $A_{0} \in
{\cal C}(R)$ is the core of sparse annular decomposition for each of
$q_{1},...,q_{P}$.
Write $S = R \setminus N(A_{0})$ and let $h \in \Mod(R)$.

We say that $h,G(q_{1}),...,G(q_{P})$ are {\em compatible} along $A_{0}$
if
\begin{enumerate}
\item $G_{0}= G_{0}(q_{1}) = ... = G_{0}(q_{P})$,
\item $h$ centralizes $G_{0}$ and,
\item $h$ is pure and pseudo-Anosov on all components of $S$.
\end{enumerate}

Our main theorem is the following.

\begin{theorem} \label{combothrm}
Suppose $h,G(q_{1}),...,G(q_{P})$ are compatible along the sparse
multi-curve $A_{0}$.
Then there exists $K_{2},...,K_{P} \geq 0$ such that for $k_{i} \geq K_{i}
+ k_{i-1}$, $i = 2,...,P$, and $k_{1} = 0$,
$$G(q_{1}) *_{G_{0}} h^{k_{2}}G(q_{2})h^{-k_{2}} *_{G_{0}} \cdots
*_{G_{0}} h^{k_{P}}G(q_{P})h^{-k_{P}} \hookrightarrow \Mod(R)$$
is an embedding.
Moreover, every element not conjugate into an elliptic or parabolic
subgroup of any $h^{k_{i}}G(q_{i})h^{-k_{i}}$ is pseudo-Anosov.
\end{theorem}

This obviously implies Theorem \ref{fakecombothrm}, where ``sufficiently
complicated'' simply means a sufficiently high power of a pure automorphism
in the centralizer of $G_{0}$ being pseudo-Anosov on all components of
$S$.
Appealing to the examples described in Section \ref{exsect}, we have the
following.\\

\noindent
{\bf Corollary \ref{sfcegrpcor}}
{\em For every closed surface $R$ of genus $g \geq 2$, there exist
subgroups of $\Mod(R)$ isomorphic to the fundamental group of a closed
surface (of genus $2g$) for which all but one conjugacy class of elements
(up to powers) is pseudo-Anosov.}\\

\noindent
{\em Proof.} Let $G(q_{1}) = G(q_{2}) = G(\xi_{g})$ from Corollary
\ref{onepunctcor}.
Since $A_{0}$ is sparse, there exists an automorphism $h \in \Mod(R)$
which is pseudo-Anosov on all component of $S$ (in fact, $S$ is connected
in all of these examples).
$G_{0}$ is generated by a multi-twist about $A_{0}$, and so $h$ certainly
centralizes $G_{0}$.
It follows that $h,G(q_{1}),G(q_{2})$ are compatible along $A_{0}$.
Note that $G_{0}$ is the only conjugacy class of parabolic subgroups in
$G(\xi_{g})$.

Theorem \ref{combothrm} therefore implies that for sufficiently large $k$
$$G(\xi_{g}) *_{G_{0}} h^{k} G(\xi_{g}) h^{-k}$$
embeds in $\Mod(R)$ and every element not conjugate into the cyclic
subgroup $G_{0}$ is pseudo-Anosov.
On the other hand
$$G(\xi_{g}) *_{G_{0}} h^{k} G(\xi_{g}) h^{-k} \cong \pi_{1}(S_{g,1})
*_{\mathbb Z} \pi_{1}(S_{g,1}) \cong \pi_{1}(S_{2g,0})$$
where the ${\mathbb Z}$ in the amalgamated product is the peripheral
subgroup of each. $\Box$\\

The proof of Theorem \ref{combothrm} is given in \S \ref{proofsect}.  This
will require some preliminary notation and results which occupy the next
few sections.

%For any sporadic surface, its mapping class group is either virtually
%free or finite.
%In particular, it does not contain the fundamental group of a closed
%surface of genus $\geq 2$.

%%%%%%%%%%%%%%%%%%%%%%%%%%%%%%%%%%%%%%%%%%%%%%%
%%%%%%%%%%%%%%%%%%%%%%%%%%%%%%%%%%%%%%%%%%%%%%%
%%%%%%%%%%%%%%%%%%%%%%%%%%%%%%%%%%%%%%%%%%%%%%%
\section{More on geodesics} \label{sadtracsect}
%%%%%%%%%%%%%%%%%%%%%%%%%%%%%%%%%%%%%%%%%%%%%%%
%%%%%%%%%%%%%%%%%%%%%%%%%%%%%%%%%%%%%%%%%%%%%%%
%%%%%%%%%%%%%%%%%%%%%%%%%%%%%%%%%%%%%%%%%%%%%%%

For the remainder of this section, fix $q \in A_{2}{\cal T}(R)$.

%%%%%%%%%%%%%%%%%%%%%%%%%%%%%%%%%%%%%%%%%%%%%%%%%%
%%%%%%%%%%%%%%%%%%%%%%%%%%%%%%%%%%%%%%%%%%%%%%%%%%
\subsection{Bad singularities} \label{badsingsect}
%%%%%%%%%%%%%%%%%%%%%%%%%%%%%%%%%%%%%%%%%%%%%%%%%%
%%%%%%%%%%%%%%%%%%%%%%%%%%%%%%%%%%%%%%%%%%%%%%%%%%

For any multi-curve $A$, we define $\sigma(A)$ to be the $q$-geodesic
representative of $A$.
This is a union of saddle connections and non-singular closed geodesics
with pairwise disjoint interiors.
We make this unique by requiring all components of $A$ which are homotopic
to nonsingular geodesics to be represented by a geodesic in the interior
of the corresponding annulus, equidistant from both boundary components.
%$\sigma(\alpha)$ is an example of a saddle track.

Since $\sigma(A)$ is a union of geodesics, if we follow an arc of any
component as it enters and exits a singularity $p$, we see that it must
make an angle at least $\pi$ on both sides, unless $p$ is a puncture
point.
In this case it must make an angle at least $\pi$ on the side opposite the
puncture.

\begin{figure}[htb]
\centerline{}
\begin{center}
\ \psfig{file=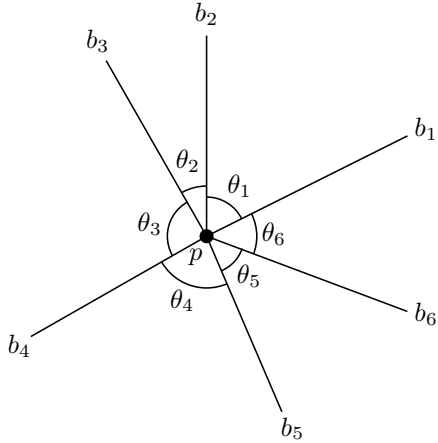,height=2truein}
\caption{Saddle connection ends} \label{saddlendsfig}
\end{center}
%%%%%%%%
  \setlength{\unitlength}{1in}
  \begin{picture}(0,0)(0,0)
    \put(4,2.07){$b_{1}$}
    \put(2.85,2.67){$b_{2}$}
    \put(2.28,2.53){$b_{3}$}
    \put(1.87,.93){$b_{4}$}
    \put(3.3,.5){$b_{5}$}
    \put(4,1.1){$b_{6}$}
    \put(3.02,1.77){$\theta_{1}$}
    \put(2.75,1.9){$\theta_{2}$}
    \put(2.55,1.57){$\theta_{3}$}
    \put(2.72,1.17){$\theta_{4}$}
    \put(3.07,1.3){$\theta_{5}$}
    \put(3.2,1.52){$\theta_{6}$}
    \put(2.82,1.42){$p$}
  \end{picture}
%%%%%%%%
\end{figure}

At each singularity $p$ of $\sigma(A)$, we consider the ends of those
saddle connections meeting $p$.
By an end we simply mean a component of a saddle connection intersected
with a small neighborhood of $p$.
We cyclically order these ends as we encounter them by encircling $p$ in a
counter-clockwise direction.
Let us write these in order as $b_{1},...,b_{r}$ (see Figure
\ref{saddlendsfig}).
Between each consecutive pair $b_{i}$ and $b_{i+1}$ (mod $r$), we have the
angle $\theta_{i}$.
We say that the singularity $p$ is {\em bad} if $\theta_{i} < \pi$ for
each $i=1,...,r$.

\vspace{.2truein}

\begin{lemma} \label{badsadlem}
Given $A \in {\cal C}(R)$, $\sigma(A)$ contains no bad singularities.
\end{lemma}

\begin{figure}[htb]
\centerline{}
\begin{center}
\ \psfig{file=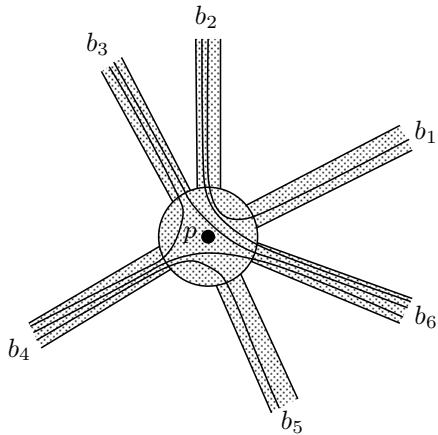,height=2truein}
\caption{Neighborhood of a singularity and its saddle connection ends}
\label{singneighfig}
\end{center}
%%%%%%%%
  \setlength{\unitlength}{1in}
  \begin{picture}(0,0)(0,0)
    \put(4,2.07){$b_{1}$}
    \put(2.85,2.67){$b_{2}$}
    \put(2.28,2.53){$b_{3}$}
    \put(1.87,.93){$b_{4}$}
    \put(3.3,.55){$b_{5}$}
    \put(4,1.1){$b_{6}$}
    \put(2.79,1.55){$p$}
  \end{picture}
%%%%%%%%
\end{figure}

\noindent
{\em Proof. } This is an outer-most arc argument.
Suppose $\sigma(A)$ contains a bad singularity.
Take a small neighborhood of $p$ and of $b_{1},...,b_{r}$ and by small
homotopy we assume that $A$ is embedded in the union of these
neighborhoods (see Figure \ref{singneighfig}).
The neighborhood of $p$ is a disk and the neighborhoods of the $b_{i}$'s
define intervals $I_{i}$ on the boundary of this disk.
The intersection of $A$ with the disk is a collection of arcs with
endpoints on the intervals.

\vspace{.2truein}

There are finitely many arcs, and so we may consider one which is
outer-most.
The endpoints of this arc are in adjacent intervals $I_{i},I_{i+1}$.
Therefore the arc lies in a component of $A$ whose geodesic representative
enters $p$ along $b_{i}$ and exits along $b_{i+1}$.
Because $\theta_{i} < \pi$, this makes an angle less than $\pi$.
Therefore $\sigma(A)$ cannot be geodesic, providing a contradiction,
unless $p$ is a puncture cut off by the arc.
In this case, we choose a different outer-most arc (or if this is the only
arc, we think of it as outer-most on the other side).
As we have already cut off the puncture with the first arc, the new
outer-most arc cannot cut it off, and we arrive at the same contradiction.
$\Box$

%%%%%%%%%%%%%%%%%%%%%%%%%%%%%%%%%%%%%%%%%%%%%%%%%%%%%%%%%%%
%%%%%%%%%%%%%%%%%%%%%%%%%%%%%%%%%%%%%%%%%%%%%%%%%%%%%%%%%%%
%\subsection{Spines and poison graphs} \label{poisonsadsect}
%%%%%%%%%%%%%%%%%%%%%%%%%%%%%%%%%%%%%%%%%%%%%%%%%%%%%%%%%%%
%%%%%%%%%%%%%%%%%%%%%%%%%%%%%%%%%%%%%%%%%%%%%%%%%%%%%%%%%%%
%Suppose that $\sigma_{0}$ and $\sigma_{1}$ are two distinct spines and
%that $\sigma_{0}'$ is a nonempty union of components of $\sigma_{0}$.
%We call $\sigma_{0}' \cup \sigma_{1}$ a {\em poison graph}.
%
%\begin{lemma} \label{2asadlem}
%Given $A \in {\cal C}(R)$, $\sigma(A)$ contains no poison graph.
%\end{lemma}
%
%\noindent
%{\em Proof.} Suppose $\sigma(A)$ contains a poison graph $\sigma_{0,1}$
%defined from $\sigma_{0}'$ and $\sigma_{1}$ as above.
%
%For any spine $\sigma$ and any singular point $p$, all the saddle
%connection ends at $p$ are evenly spaced at angles $\pi$.
%Since $\sigma_{0}$ and $\sigma_{1}$ are distinct spines, their saddle
%connection ends are perfectly interlaced (see Figure \ref{2asaddlefig}).
%In particular, any singular point of $\sigma_{0,1}$ which $\sigma_{0}'$
%meets is bad.
%Therefore, $\sigma(A)$ has a bad singularity, which by Lemma
%\ref{badsadlem} is a contradiction. $\Box$
%
%\begin{figure}[htb]
%\centerline{}
%\begin{center}
%\ \psfig{file=new2asaddle.eps,height=2.5truein}
%\caption{Local picture near the singularity of a poison graph}
%\label{2asaddlefig}
%\end{center}
%\end{figure}
%%%%%%%%%%%%%%%%%%%%%%%%%%%%%%%%%%%%%%%%%%%%%%%%%%%%%%%%%%%%%%%
%%%%%%%%%%%%%%%%%%%%%%%%%%%%%%%%%%%%%%%%%%%%%%%%%%%%%%%%%%%%%%%
\subsection{Subsurface projections and spines} \label{asadsect}
%%%%%%%%%%%%%%%%%%%%%%%%%%%%%%%%%%%%%%%%%%%%%%%%%%%%%%%%%%%%%%%
%%%%%%%%%%%%%%%%%%%%%%%%%%%%%%%%%%%%%%%%%%%%%%%%%%%%%%%%%%%%%%%

Given any annular decomposition of $R$ for $q$ with core $A_{0}$, let
$\sigma_{0}$ denote the union of singular leaves parallel to $A_{0}$.
Equivalently, this is the union of the boundary components of maximal
annuli in the annular decomposition.
Call $\sigma_{0}$ a {\em spine} parallel to $A_{0}$.

Fix a spine $\sigma_{0}$ parallel to a sparse core $A_{0}$ and let
$\sigma_{0,1},..., \sigma_{0,N}$ denote the components of $\sigma_{0}$.
For each $j = 1,..., N$, define
$$\Theta(q)_{j} = \{ A \in {\cal C}(R) \setminus Z'(A_{0}) \, | \,
\sigma(A) \not \supset \sigma_{0,j} \}$$
and set
$$\Theta(q) = \bigcap_{j=1}^{N} \Theta(q)_{j}$$

\begin{figure}[htb]
\centerline{}
\begin{center}
\ \psfig{file=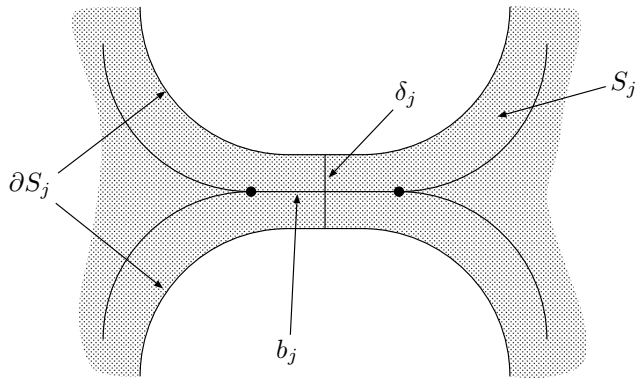,height=2truein}
\caption{The arc $\delta_{j}$} \label{globalboundfig}
\end{center}
%%%%%%%%
  \setlength{\unitlength}{1in}
  \begin{picture}(0,0)(0,0)
    \put(4.45,2.15){$S_{j}$}
    \put(2.7,.75){$b_{j}$}
    \put(3.32,2.1){$\delta_{j}$}
    \put(1.3,1.62){$\partial S_{j}$}
  \end{picture}
%%%%%%%%
\end{figure}

Let $S = R \setminus N(A_{0})$ with components $S_{1},..., S_{N}$.
We may order these components so that $\sigma_{0,j} \subset S_{j}$.
For each $j = 1,...,N$, take $\delta_{j}$ to be a geodesic arc in $S_{j}$
with endpoints on $\partial S_{j}$ and intersecting the interior of a
single saddle connection, $b_{j}$, of $\sigma_{0,j}$ (see Figure
\ref{globalboundfig}).

\begin{proposition} \label{projboundprop}
For each $j = 1,...,N$,
$$\pi_{S_{j}}(\Theta(q)) \subset \pi_{S_{j}}(\Theta(q)_{j}) \subset
B(\delta_{j},2)$$
\end{proposition}

\noindent
{\em Proof.}  The first inclusion is obvious, so all we need prove is the
second.
Before we begin, we note that for any sufficiently small $\epsilon > 0$,
we can take $S_{j}$ to be the $\epsilon$-neighborhood of $\sigma_{0,j}$.

There is a linear deformation retraction $H_{t}$, $t \in [0,1]$, of
$S_{j}$ onto $\sigma_{0,j}$ which slides points along a perpendicular
geodesic arc from $\partial S_{j}$ to $\sigma_{0,j}$.
We assume that $H_{t}$ is an embedding for $t \in [0,1)$.
Moreover, $\delta_{j}$ may be chosen to be $H_{1}^{-1}(p)$ for any $p \in
\intr(b_{j})$ (the interior of $b_{j}$), and this is isotopic to any other
choice of $\delta_{j}$.

We fix some $A \in \Theta(q)_{j}$.
By definition, $\sigma_{0,j} \not \subset \sigma(A)$.
Suppose first that $b_{j} \not \subset \sigma(A)$.
We claim that for an appropriate choice of $p \in \intr(b_{j})$,
$\delta_{j} = H_{1}^{-1}(p)$ is disjoint from $\sigma(A)$.

To prove the claim, we note that because $b_{j} \not \subset \sigma(A)$,
we can find $p \in \intr(b_{j})$ and $\epsilon' > 0$ such that
$$d_{q}(p,\sigma(A)) > 2 \epsilon'$$
Now, take any $0< \epsilon \leq \epsilon'$ which defines $S_{j}$ as the
$\epsilon$-neighborhood of $\sigma_{0,j}$.
Let $\gamma$ be an arc of $b_{j}$ centered at $p$ of radius no more than
$\epsilon$.
Then $H_{1}^{-1}(\gamma)$ contains $\delta_{j}$ and is contained in the $2
\epsilon'$-neighborhood of $p$.
In particular, this set (and hence $\delta_{j}$) is disjoint from
$\sigma(A)$ as claimed.

Thus, if $b_{j} \not \subset \sigma_{0,j}$, then $\delta_{j}$ is disjoint
from $\sigma(A)$ and in particular
$$d(\delta_{j},\pi_{S_{j}}(A)) \leq 1$$
for any allowable choice of projection $\pi_{S_{j}}(A)$.

So to prove the proposition, we assume $b_{j} \subset \sigma(A)$.
Then by hypothesis there is another branch $b_{j}'$ of $\sigma_{0,j}$ such
that $b_{j}' \not \subset \sigma(A)$.
Let $\delta_{j}'$ be an analogous arc for $b_{j}'$.
The previous argument implies that $d(\delta_{j}',\pi_{S_{j}}(A)) \leq 1$
for any allowable choice of projection.
$\delta_{j}$ and $\delta_{j}'$ are disjoint however, and so
$d(\delta_{j},\pi_{S_{j}}(A)) \leq 2$. $\Box$

%%%%%%%%%%%%%%%%%%%%%%%%%%%%%%%%%%%%%%%%%%%%%%%%%%%%%%%%%%%%
%%%%%%%%%%%%%%%%%%%%%%%%%%%%%%%%%%%%%%%%%%%%%%%%%%%%%%%%%%%%
\subsection{$G(q)$ action on geodesics} \label{gqactsadsect}
%%%%%%%%%%%%%%%%%%%%%%%%%%%%%%%%%%%%%%%%%%%%%%%%%%%%%%%%%%%%
%%%%%%%%%%%%%%%%%%%%%%%%%%%%%%%%%%%%%%%%%%%%%%%%%%%%%%%%%%%%

$G(q)$ acts on the set of saddle connections and non-singular closed
geodesics since it is acting by affine homeomorphisms.
Taking geodesic representatives is obviously equivariant with respect to
this action and the action on ${\cal C}(R)$.
That is
\begin{equation} \label{sadacteqn}
\sigma(\phi(A)) = \phi(\sigma(A))
\end{equation}
for any $A \in {\cal C}(R)$ and any $\phi \in G(q)$.

This implies
\begin{lemma} \label{thetainvtlem}
$\Theta(q)$ is invariant under $G_{0}(q)$.
\end{lemma}

\noindent
{{\em Proof.} $\sigma_{0}$ is invariant under $G_{0}(q)$ although some
elements may permute the components.
Let $A \in {\cal C}(R)$ and $\phi \in G_{0}(q)$ and let
$\phi(\sigma_{0,j}) = \sigma_{0,j'}$.
Then we have
$$\sigma_{0,j'} = \phi(\sigma_{0,j}) \subset \phi(\sigma(A) =
\sigma(\phi(A))  \Leftrightarrow \sigma_{0,j} \subset \sigma(A)$$
Since $A_{0}$ is also invariant under $G_{0}(q)$, so is $Z'(A_{0})$.
The lemma now follows from the definition of $\Theta(q)$. $\Box$

%%%%%%%%%%%%%%%%%%%%%%%%%%%%%%%%%%%%%%%%
%%%%%%%%%%%%%%%%%%%%%%%%%%%%%%%%%%%%%%%%
%%%%%%%%%%%%%%%%%%%%%%%%%%%%%%%%%%%%%%%%
\section{The ``horoball'' sets} \label{hballsect}
%%%%%%%%%%%%%%%%%%%%%%%%%%%%%%%%%%%%%%%%
%%%%%%%%%%%%%%%%%%%%%%%%%%%%%%%%%%%%%%%%
%%%%%%%%%%%%%%%%%%%%%%%%%%%%%%%%%%%%%%%%

We wish to define the analog of the horoballs used in the Kleinian group
example of \S \ref{kleinsect}.
These sets will depend on $q$ and a sparse core $A_{0}$ of an
annular decomposition for $q$.

We define $H(q) = H(q,A_0)$ by the equation
$$H(q) = \bigcap_{j=1}^{N} \pi_{S_{j}}^{-1}\left( {\cal A}^{0}(S_{j})
\setminus B(\delta_{j},2)  \right) \subset {\cal C}(R) \setminus Z'(A_{0})$$
where $\delta_{j}$ is as in \S \ref{asadsect} and Proposition
\ref{projboundprop}.

\begin{proposition} \label{horobehavprop}
For every $j = 1,...,N$, we have
$$H(q) \cap \Theta(q)_{j} = \emptyset$$
(in particular, $H(q) \cap \Theta(q) = \emptyset$), and for every $\phi
\in G(q) \setminus G_{0}$
$$\phi(H(q)) \subset \Theta(q)$$
\end{proposition}

\noindent
{\em Proof. } By Proposition \ref{projboundprop}, for every $j = 1,...,N$,
we have
$$\left( {\cal A}^{0}(S_{j}) \setminus B(\delta_{j},2) \right) \cap
\pi_{S_{j}}(\Theta(q)_{j}) = \emptyset$$
Therefore
$$\begin{array}{rcl}
H(q) \cap \Theta(q)_{j} & = & \left( \bigcap_{k=1}^{N} \pi_{S_{k}}^{-1}
\left( {\cal A}^{0}(S_{k}) \setminus B(\delta_{k},2)  \right) \right) \cap
\Theta(q)_{j} \\
 & \subset & \pi_{S_{j}}^{-1} \left( {\cal A}^{0}(S_{j}) \setminus
B(\delta_{j},2)  \right) \cap \pi_{S_{j}}^{-1}
\left(\pi_{S_{j}}(\Theta(q)_{j})\right) \\
 & = & \pi_{S_{j}}^{-1}\left(  \left( {\cal A}^{0}(S_{j}) \setminus
B(\delta_{j},2)  \right) \cap \pi_{S_{j}}(\Theta(q)_{j})\right)\\
 & = & \emptyset \\ \end{array}$$
which proves the first assertion.

To prove the second, we suppose that $A \in H(q)$ and $\phi \in G(q)
\setminus G_{0}$ and show that $\phi(A) \in \Theta(q)$.
By the first part, we know $A \not \in \Theta(q)_{j}$ for any $j =
1,...,N$ and so $\sigma_{0} \subset \sigma(A)$.
Since
$$\phi(\sigma(A)) = \sigma(\phi(A))$$
by (\ref{sadacteqn}), we see that $\phi(\sigma_{0}) \subset
\sigma(\phi(A))$.
Since $\sigma_{0}$ is a spine, so is $\phi(\sigma_{0})$.
Now, suppose $\sigma_{0,j} \subset \sigma(\phi(A))$ for some $j$.
$\sigma_{0,j}$ is a component of the spine $\sigma_{0}$ and so
$\sigma(\phi(A))$ contains the $\sigma_{0,j} \cup
\phi(\sigma_{0})$. This has a bad vertex (see Figure \ref{2asaddlefig})
and so contradicts Lemma \ref{badsadlem}. Hence $\sigma_{0,j} \not \subset
\sigma(\phi(A))$ for any $j = 1,...,N$.

\begin{figure}[htb]
\centerline{}
\begin{center}
\ \psfig{file=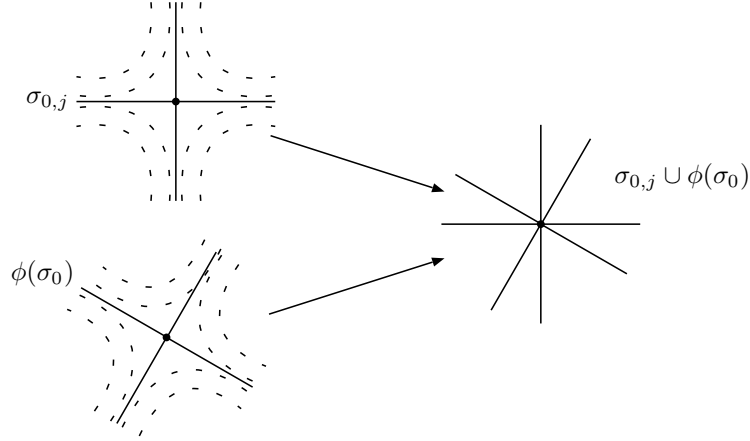,height=2.5truein}
\caption{Local picture near the singularity of 
                             $\sigma_{0,j} \cup \phi(\sigma_{0})$}
\label{2asaddlefig}
\end{center}
%%%%%%%%
  \setlength{\unitlength}{1in}
  \begin{picture}(0,0)(0,0)
    \put(4.45,2.15){$\sigma_{0,j} \cup \phi(\sigma_{0})$}
    \put(1.37,2.57){$\sigma_{0,j}$}
    \put(1.29,1.64){$\phi(\sigma_{0})$}
  \end{picture}
%%%%%%%%
\end{figure}

To prove that $\phi(A) \in \Theta(q)$, all that remains is to show that
$\phi(A) \not \in Z'(A_{0})$.
Again, we note that $\phi(\sigma_{0})$ is a spine and it therefore
(nontrivially and transversely) intersects any nonsingular geodesic which
its associated core, $\phi(A_{0})$, intersects.
However, $\phi(A_{0})$ and $A_{0}$ are distinct cores and so Lemma
\ref{allinterlem} guarantees $\phi(A_{0})$ intersects every component of
$A_{0}$.
Therefore, $\phi(A)$ intersects every component of $A_{0}$ and thus
$\phi(A) \not \in Z'(A_{0})$. $\Box$

%%%%%%%%%%%%%%%%%%%%%%%%%%%%%%%%%%%%%%%%%%%%%
%%%%%%%%%%%%%%%%%%%%%%%%%%%%%%%%%%%%%%%%%%%%%
\subsection{Separating projections} \label{spaceoutsect}
%%%%%%%%%%%%%%%%%%%%%%%%%%%%%%%%%%%%%%%%%%%%%
%%%%%%%%%%%%%%%%%%%%%%%%%%%%%%%%%%%%%%%%%%%%%

We assume the hypotheses on $h,G(q_{1}),...,G(q_{P})$ and $A_{0}$ from the
statement of Theorem \ref{combothrm} in what follows.
Let $\Theta(q_{i})_{j}$ for $j = 1,...,N$ and $i = 1,...,P$ be defined as
in \S \ref{asadsect} with respect to our given multi-curve $A_{0} \in
{\cal C}(R)$.
Let $h_{j} \in \Mod(S_{j})$ for each $j = 1,...,N$ be the
$S_{j}$-component of $h$, which is pseudo-Anosov by assumption.
As a convenience, we use the flexibility in the definition of
$\pi_{S_{j}}$ (see \S \ref{subprosect}) to guarantee that it is
equivariant with respect to $\langle h \rangle$.
That is, for each $j = 1,...,N$, each $k \in {\mathbb Z}$, and each $A \in
{\cal C}(R) \setminus Z'(A_{0})$
$$\pi_{S_{j}}(h^{k}(A)) = h_{j}^{k}(\pi_{S_{j}}(A))$$
This is accomplished by defining $\pi_{S_{j}}$ on $\langle h
\rangle$-orbit representatives, then extending so that it is equivariant
(it is easy to check that this is an allowable definition for
$\pi_{S_{j}}$).

In the statement of Theorem \ref{combothrm} we require the existence of
certain numbers $K_{2},...,K_{P}$.
The following lemma will provide us with these numbers.

\begin{lemma} \label{spacedlem}
Suppose $h,G(q_{1}), ...,G(q_{P})$ are compatible along the sparse
multi-curve $A_{0}$.
Then there exists $K_{2},...,K_{P} \geq 0$ such that for $k_{i} \geq K_{i}
+ k_{i-1}$, $i = 2,...,P$, and $k_{1} = 0$
$$d(h_{j}^{k_{i}}(\pi_{S_{j}}(\Theta(q_{i}))),h_{j}^{k_{i'}}(\pi_{S_{j}}(\Theta(q_{i'}))))
\geq 5$$
for every $i,i' \in \{1,...,P \}$ with $i \neq i'$.
\end{lemma}

\noindent
{\em Proof.} Proposition \ref{projboundprop} implies that for each $j =
1,...,N$ and each $i = 1,...,P$, $\pi_{S_{j}}(\Theta(q_{i}))$ has diameter
no more than $4$.
In particular, by Theorem \ref{mmpathrm} and the triangle inequality, for
each $i = 2,...,P$ there are only finitely many integers $k_{i}$ for which
$$d(h_{j}^{k_{i}}(\pi_{S_{j}}(\Theta(q_{i}))),\pi_{S_{j}}(\Theta(q_{i'})))
< 5$$
for some $i' \neq i$.
Therefore, for each $i = 2,...,P$, there exists $K_{i} \geq 0$ such that
for every $k_{i} \geq K_{i}$ and every $i' \neq i$,
\begin{equation} \label{ispacedeqn}
d(h_{j}^{k_{i}}(\pi_{S_{j}}(\Theta(q_{i}))),\pi_{S_{j}}(\Theta(q_{i'})))
\geq 5
\end{equation}

Now set $k_{1} = 0$ and let $k_{i} \geq K_{i} + k_{i-1}$.
Suppose $i > i'$.
By induction on $i-i'$, we see that
$$k_{i} \geq \left( \sum_{l=i'+1}^{i} K_{l} \right) + k_{i'}$$
In particular, $k_{i} - k_{i'} \geq K_{i}$.

Since $h_{j}^{-k_{i'}}$ acts by isometries on ${\cal A}^{0}(S_{j})$, we
have
$$d(h_{j}^{k_{i}}(\pi_{S_{j}}(\Theta(q_{i}))),h_{j}^{k_{i'}}(\pi_{S_{j}}(\Theta(q_{i'}))))
=
d(h_{j}^{k_{i}-k_{i'}}(\pi_{S_{j}}(\Theta(q_{i}))),\pi_{S_{j}}(\Theta(q_{i'})))
\geq 5$$
by (\ref{ispacedeqn}), as required. $\Box$\\

%%%%%%%%%%%%%%%%%%%%%%%%%%%%%%%%%%%%%%%%%%%%%%%%%%%%%%%%%%%%
%%%%%%%%%%%%%%%%%%%%%%%%%%%%%%%%%%%%%%%%%%%%%%%%%%%%%%%%%%%%
%%%%%%%%%%%%%%%%%%%%%%%%%%%%%%%%%%%%%%%%%%%%%%%%%%%%%%%%%%%%
\section{Proof of Theorem \ref{combothrm}} \label{proofsect}
%%%%%%%%%%%%%%%%%%%%%%%%%%%%%%%%%%%%%%%%%%%%%%%%%%%%%%%%%%%%
%%%%%%%%%%%%%%%%%%%%%%%%%%%%%%%%%%%%%%%%%%%%%%%%%%%%%%%%%%%%
%%%%%%%%%%%%%%%%%%%%%%%%%%%%%%%%%%%%%%%%%%%%%%%%%%%%%%%%%%%%

We can now prove\\

\noindent {\bf Theorem \ref{combothrm}}
{\em Suppose $h,G(q_{1}),...,G(q_{P})$ are compatible along the sparse
multi-curve $A_{0}$.
Then there exists $K_{2},...,K_{P} \geq 0$ such that for $k_{i} \geq K_{i}
+ k_{i-1}$, $i = 2,...,P$, and $k_{1} = 0$,
$$G(q_{1}) *_{G_{0}} h^{k_{2}}G(q_{2})h^{-k_{2}} *_{G_{0}} \cdots
*_{G_{0}} h^{k_{P}}G(q_{P})h^{-k_{P}} \hookrightarrow \Mod(R)$$
is an embedding.
Moreover, every element not conjugate into an elliptic or parabolic
subgroup of any $h^{k_{i}}G(q_{i})h^{-k_{i}}$ is pseudo-Anosov.}\\

\noindent
{\em Proof.} Let $K_{2},...,K_{P}$ be as in Lemma \ref{spacedlem}, $k_{1}
= 0$, and $k_{i} \geq K_{i} + k_{i-1}$ for each $i = 2,...,P$.
To simplify the notation, we replace each $G(q_{i})$ by its conjugate
$h^{k_{i}} G(q_{i}) h^{-k_{i}}$.
By the $\langle h \rangle$-equivariance of $\pi_{S_{j}}$ (see \S
\ref{spaceoutsect}), this has the effect of replacing each
$\pi_{S_{j}}(\Theta(q_{i}))$ by
$$\pi_{S_{j}}(h^{k_{i}}(\Theta(q_{i}))) =
h_{j}^{k_{i}}(\pi_{S_{j}}(\Theta(q_{i})))$$
Therefore, by Lemma \ref{spacedlem}, we have
\begin{equation} \label{spaced2eqn}
d(\pi_{S_{j}}(\Theta(q_{i})),\pi_{S_{j}}(\Theta(q_{i'}))) \geq 5
\end{equation}

We first claim that $(\Theta(q_{1}),...,\Theta(q_{P}))$ is a proper
interactive $P$-tuple (see \S \ref{pipsect}).
Clearly each $\Theta(q_{i})$ is nonempty and Lemma \ref{thetainvtlem}
implies each of these sets is invariant under $G_{0}$.
By (\ref{spaced2eqn}) we see that for each $i \neq i'$ and each $j =
1,...,N$
$$\pi_{S_{j}}(\Theta(q_{i'})) \subset {\cal A}^{0}(S_{j}) \setminus
B(\delta_{i,j},2)$$
where $\delta_{i,j}$ is defined for $S_{j}$ in terms of $q_{i}$ as in \S
\ref{asadsect} and Proposition \ref{projboundprop} so that
\begin{equation} \label{thetainhoroeqn}
\Theta(q_{i'}) \subset \bigcap_{j=1}^{N}
\pi_{S_{j}}^{-1}(\pi_{S_{j}}(\Theta(q_{i'}))) \subset H(q_{i})
\end{equation}
Thus, by the first part of Proposition \ref{horobehavprop} (or by
(\ref{spaced2eqn}) directly), we see that
$$\Theta(q_{i}) \cap \Theta(q_{i'}) = \emptyset$$
for each $i \neq i'$.

Now, the second part of Proposition \ref{horobehavprop} combined with
(\ref{thetainhoroeqn}) shows that for every $i$ and each $\phi \in
G(q_{i}) \setminus G_{0}$
\begin{equation} \label{pingeqn}
\phi(\Theta(q_{i'})) \subset \phi(H(q_{i})) \subset \Theta(q_{i})
\end{equation}
for every $i' \neq i$.
Therefore $\Theta(q_{1}),...,\Theta(q_{P})$ satisfies properties 1--4 in
the definition of proper interactive $P$-tuple.

To prove that this also satisfies property 5, for every $i$ choose any
$\phi_{i} \in G(q_{i}) \setminus G_{0}$.
We claim that
\begin{equation} \label{wheresalphaeqn}
\phi_{i}(A_{0}) \in \Theta(q_{i})
\end{equation}
This is because $\phi_{i}(A_{0})$ is a core of an annular decomposition
different than $A_{0}$, so $\sigma(\phi_{i}(A_{0}))$ (with respect to
$q_{i}$) contains no saddle connections and therefore no component of
$\sigma_{0}$.
Moreover, by Lemma \ref{allinterlem}, $\phi_{i}(A_{0}) \not \in
Z'(A_{0})$.
Define $\theta_{i} = \phi_{i}(A_{0})$.

For every $i' \neq i$ and any $\tilde{\phi} \in G(q_{i}) \setminus G_{0}$,
we have
$$\theta_{i} \not \in \tilde{\phi}(\Theta(q_{i'}))$$
For, as we saw in the proof of Proposition \ref{horobehavprop}, given any
$A \in \tilde{\phi}(\Theta(q_{i'}))$, $\sigma(A)$ contains a spine (with
respect to $q_{i}$), namely, $\tilde{\phi}(\sigma_{0})$.
However, $\theta_{i}$ contains no spine.
We have thus proved that $\Theta(q_{1}),...,\Theta(q_{P})$ also satisfies
property 5 and hence is a proper interactive $P$-tuple.
Proposition \ref{pipprop} implies the required injectivity.\\

We now prove the statement about pseudo-Anosov elements.
By Theorem \ref{affactthrm} we know that in each factor, any element not
in an elliptic or parabolic subgroup is pseudo-Anosov.
It suffices therefore to consider an element $\phi$ which is not conjugate
into any factor.
We will assume that $\phi$ is reducible and arrive at a contradiction.
By the structure theory of amalgamated free products, since $\phi$ 
is not conjugate into a factor, 
$\phi$ must have infinite order.
Theorem \ref{autclassthrm} will then imply that $\phi$ is pseudo-Anosov.

After conjugating, we may write
$$\phi = \phi_{i_{1}} \phi_{i_{2}} \cdots \phi_{i_{r}}$$
with $\phi_{i_{k}} \in G(q_{i_{k}}) \setminus G_{0}$ and $i_{k} \neq
i_{k+1}$ (mod $r$) for each $k = 1,...,r$.
By passing to an appropriate power if necessary, we may also assume that
$\phi$ is pure.
In what follows, we let $A$, $\lambda_{s,1}$, ... , $\lambda_{s,M}$,
$\lambda_{u,1}$, ... , $\lambda_{u,M}$, $\psi_{s}$, $\psi_{u}$,
$\Psi_{s}$, and $\Psi_{u}$ be as in \S \ref{identpasect}.

We consider the image of $A_{0}$ under $\phi$ and its iterates.
First, by (\ref{wheresalphaeqn}), $\phi_{i_{r}}(A_{0}) \in
\Theta(q_{i_{r}})$.
It follows from (\ref{pingeqn}) that $\phi_{i_{r-1}}(\phi_{i_{r}}(A_{0}))
\in \Theta(q_{i_{r-1}})$.
Applying the rest of the word for $\phi$, by repeatedly appealing to
(\ref{pingeqn}) and induction, we see
$$\phi(A_{0}) \in \Theta(q_{i_{1}})$$
Similarly,
$$\phi^{-1}(A_{0}) \in \Theta(q_{i_{r}})$$

Iterating $\phi$ and $\phi^{-1}$ and applying (\ref{pingeqn}) and
induction we find that
\begin{equation} \label{whereweliveeqn}
\phi^{m}(A_{0}) \in \Theta(q_{i_{1}}) \mbox{ and } \phi^{-m}(A_{0}) \in
\Theta(q_{i_{r}})
\end{equation}
for every $m \geq 1$.

So, combining (\ref{spaced2eqn}) and (\ref{whereweliveeqn}) we see that
for every $j = 1,...,N$ and $m \geq 1$
\begin{equation} \label{faraparteqn}
d(\pi_{S_{j}}(\phi^{m}(A_{0})),\pi_{S_{j}}(\phi^{-m}(A_{0}))) \geq 5
\end{equation}

Now consider the limits of these multi-curves in $\prj\ML_{0}(R)$
$$[\lambda_{0}^{s}] = \lim_{m \rightarrow \infty} \phi^{m}(A_{0}) \mbox{
and } [\lambda_{0}^{u}] = \lim_{m \rightarrow \infty} \phi^{-m}(A_{0})$$
($A_{0}$ is assigned the projective class of transverse counting measure,
although any transverse measure of full support will do).
Theorem \ref{autotypethrm} implies that if $A_{0}$ nontrivially intersects
$\psi_{s}$, then $[\lambda_{0}^{s}] \in \prj\Psi_{s}(\phi)$ and
$[\lambda_{0}^{u}] \in \prj\Psi_{u}(\phi)$ (because $A_{0}$ is a
multi-curve it intersects $\psi_{s}$ if and only if it intersects
$\psi_{u}$).
However, because $\phi(A_{0}) \in \Theta(q_{i_{1}})$ and $A_{0} \not \in
\Theta(q_{i_{1}})$, the comment following Theorem \ref{autotypethrm}
implies that $A_{0}$ nontrivially intersects $\psi_{s}$.
In particular, we find $[\lambda_{0}^{s}], [\lambda_{0}^{u}] \in
\prj\ML_{0}(R) \setminus Z(\partial S_{j})$ for some $j$.

By Proposition \ref{coarsecontprop}, the triangle inequality, and
(\ref{faraparteqn}), we see that for sufficiently large $m$,
$$\begin{array} {rcll}
d(\pi_{S_{j}}([\lambda_{0}^{s}]),\pi_{S_{j}}([\lambda_{0}^{u}])) & \geq &
d(\pi_{S_{j}}(\phi^{m}(A_{0})),\pi_{S_{j}}(\phi^{-m}(A_{0}))) & \\
 &      & - \,
d(\pi_{S_{j}}(\phi^{m}(A_{0})),\pi_{S_{j}}([\lambda_{0}^{s}])) -
d(\pi_{S_{j}}([\lambda_{0}^{u}]),\pi_{S_{j}}(\phi^{-m}(A_{0}))) \\
 & \geq  & 5 - 2 \quad = \quad 3 & \\ \end{array}$$
Therefore $\lambda_{0}^{s},\lambda_{0}^{u}$ bind when restricted to
$S_{j}$.
Said differently, the intersections of $|\lambda_{0}^{s}|$ and
$|\lambda_{0}^{u}|$ with $S_{j}$ nontrivially intersect every essential
simple closed curve and essential arc in $S_{j}$.

Suppose $|\lambda_{0}^{s}|$ misses some component of $A_{0}$.
The connectivity of $R$ implies that, among all components of $A_{0}$
which $|\lambda_{0}^{s}|$ fails to intersect, there is (at least) one which
is the boundary of a component $S_{j}$ of $S$ for which $\lambda_{0}^{s},\lambda_{0}^{u}$
bind when restricted (as in the previous paragraph).
Call this component $a$.\\

\noindent
{\bf Claim.} $a$ is also a component of $A$.\\

\noindent
{\em Proof of claim.} Because $\lambda_{0}^{s},\lambda_{0}^{u}$ bind when
restricted to $S_{j}$, and since they miss $a$, the boundary of the
complementary region of $|\lambda_{0}^{s}| \cup |\lambda_{0}^{u}|$ must
have a component parallel to $a$ (see Figure \ref{boundaryparallelfig}).

\begin{figure}[htb]
\centerline{}
\centerline{}
\begin{center}
\ \psfig{file=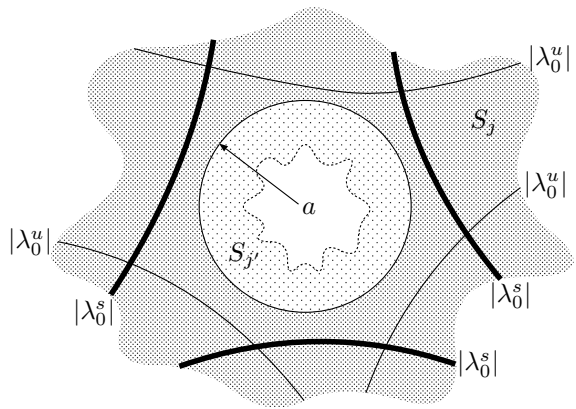,height=2.25truein}
\caption{$|\lambda_{0}^{s}| \cup |\lambda_{0}^{u}|$ near $a$.}
\label{boundaryparallelfig}
\end{center}
%%%%%%%%
  \setlength{\unitlength}{1in}
  \begin{picture}(0,0)(0,0)
    \put(3.8,2.15){$S_{j}$}
    \put(2.53,1.45){$S_{j'}$}
    \put(2.92,1.7){$a$}
    \put(3.73,.88){$|\lambda_{0}^{s}|$}
    \put(1.72,1.15){$|\lambda_{0}^{s}|$}
    \put(3.9,1.24){$|\lambda_{0}^{s}|$}
    \put(1.38,1.52){$|\lambda_{0}^{u}|$}
    \put(4.08,1.82){$|\lambda_{0}^{u}|$}
    \put(4.1,2.48){$|\lambda_{0}^{u}|$}
  \end{picture}
%%%%%%%%
\end{figure}

The components of $|\lambda_{0}^{s}| \cup |\lambda_{0}^{u}|$ are
components of $\psi_{s} \cup \psi_{u}$.
Any such component is either a component of $A$ or a set of the form
$|\lambda_{s,j}| \cup |\lambda_{u,j}|$ for $j = 1,...,M$.
So, the component of the boundary of the complementary region above is a
component of the boundary of the complementary region for some
$|\lambda_{s,j}| \cup |\lambda_{u,j}|$.
However, all such boundary components which are essential in $R$ are
parallel to the boundaries of subsurfaces on which $\phi$ is
pseudo-Anosov.
Therefore, $a$ is such a boundary, and hence a component of $A$ (the
boundary of a component on which $\phi$ is pseudo-Anosov is necessarily in
the canonical reducing system).
This proves the claim.\\

Since $\phi(A_{0})$ transversely intersects {\em every} component of
$A_{0}$ (it lies in the complement of $Z'(A_{0})$), it follows that
$\phi(a) \neq a$.
However, this contradicts the fact that $a$ is a reducing curve and was
thus assumed to be fixed (by the purity assumption).
This proves that $\phi$ is pseudo-Anosov and completes the proof the
theorem. $\Box$

%%%%%%%%%%%%%%%%%%%%%%%%%%%%%%%%%%%%%%%%%%%%%%%%%%%%%%
%%%%%%%%%%%%%%%%%%%%%%%%%%%%%%%%%%%%%%%%%%%%%%%%%%%%%%
\subsection{Other combinations} \label{othercombosect}
%%%%%%%%%%%%%%%%%%%%%%%%%%%%%%%%%%%%%%%%%%%%%%%%%%%%%%
%%%%%%%%%%%%%%%%%%%%%%%%%%%%%%%%%%%%%%%%%%%%%%%%%%%%%%

The arguments in the proof of Theorem \ref{combothrm} carry over verbatim
to more general settings.
For example, one can form the amalgamated product of countably many groups
$\{ G(q_{i}) \}_{i = 1}^{\infty}$ conjugated by (appropriate powers of)
distinct pure automorphisms $\{ h(i) \}_{i=1}^{\infty}$ each centralizing
$G_{0}$.
However, the statement seems sufficiently messy that we have only stated
it in the finite case.

One example of the countable case which seems interesting enough to
mention explicitly is when $G(q_{i}) = G(q)$ for some fixed $q$.
In this case, we have

\begin{theorem} \label{countcombothrm}
Suppose $h$ and $G(q)$ are compatible along $A_{0}$.
Then there exists $k > 0$ such that
$$G = \cdots h^{-2k}G(q)h^{2k} *_{G_{0}} h^{-k}G(q)h^{k}  *_{G_{0}} G(q)
*_{G_{0}} h^{k}G(q)h^{-k} *_{G_{0}} h^{2k}G(q)h^{-2k} \cdots
\hookrightarrow \Mod(R)$$
is an embedding and any element not conjugate into an elliptic or
parabolic subgroup of $G(q)$ is pseudo-Anosov.
\end{theorem}

\noindent
{\em Proof.} All that is needed is the analog of Lemma \ref{spacedlem} to
show that there exists $k > 0$ so that
$$d(h^{rk}(\pi_{S_{j}}(\Theta(q))),h^{r'k}(\pi_{S_{j}}(\Theta(q)))) \geq
5$$
for all $r, r' \in {\mathbb Z}$ and $r \neq r'$.
This follows from Proposition \ref{projboundprop}, Theorem \ref{mmpathrm},
and the fact that $h_{j}$ acts on ${\cal A}^{0}(S_{j})$ by isometries.
The rest of the proof is identical to the proof of Theorem
\ref{combothrm}.$\Box$\\

We now consider the trivial HNN extension
$$G(q) *_{G_{0}}$$
(we are extending by ${\mathbb Z}$ over $G_{0}$ and the stable letter,
$\tau$, acts trivially on $G_{0}$).
We write
$$\eta:G(q) *_{G_{0}} \rightarrow {\mathbb Z}$$
for the epimorphism obtained by sending $G(q)$ to $0$ and $\tau$ to $1$.

There is a homomorphism
$$\delta: G(q) *_{G_{0}} \rightarrow \Mod(R)$$
extending the inclusion of $G(q)$ obtained by sending $\tau$ to $h^{k}$.
Then the kernel of $\eta$ is isomorphic to the group $G$ from Theorem
\ref{countcombothrm} in such a way that the restriction of $\delta$ is the
inclusion of $G$ into $\Mod(R)$.

\begin{corollary}
With the notation as above,
$$\delta: G(q) *_{G_{0}} \hookrightarrow \Mod(R)$$
is an embedding.
Moreover, any element not conjugate into an elliptic or parabolic subgroup
or into $\langle h^{k} \rangle$ is pseudo-Anosov.
\end{corollary}

\noindent
{\em Proof.} We have the short exact sequence
$$1 \rightarrow G \rightarrow G(q)*_{G_{0}} \rightarrow {\mathbb Z}
\rightarrow 1$$
which splits by the homomorphism sending $1$ to $\tau$.
By construction, the splitting homomorphism composed with $\delta$ is
injective with image $\langle h^{k} \rangle$.

Any element in $\phi \in G(q)*_{G_{0}} \setminus \, \, \{ 1 \}$ is a
product $\tau^{n} \gamma$ where $\gamma \in G$ and $n \in {\mathbb Z}$.
If $\delta(\phi) = 1$ then $n \neq 0$ and $\gamma \neq 1$ since $\delta$
restricted to each of $G$ and $\langle \tau \rangle$ is injective.
Moreover
$$\gamma = \delta(\tau^{-n}\tau^{n}\gamma) = \delta(\tau^{-n}\phi) =
\delta(\tau^{-n}) = h^{-nk}$$
This is impossible since the element $\gamma \in G$ is either
pseudo-Anosov, a root of a positive multi-twist, or of finite order by
Theorem \ref{countcombothrm} and Theorem \ref{affactthrm}, whereas
$h^{-nk}$ is reducible with all pseudo-Anosov components by hypothesis.

The proof of the pseudo-Anosov statement follows by an argument similar to
that in the proof of Theorem \ref{combothrm} where we conjugate our
element to have the form
$$\phi = h^{nk} \phi_{i_{1}} \cdots \phi_{i_{r}}$$
with $\phi_{i_{j}} \in h^{i_{j}k}G(q)h^{-i_{j}k} \setminus G_{0}$, $i_{j}
\neq i_{j+1}$ (mod $r$), and $r$ minimal.
For then, $\phi^{m}(A_{0}) \in h^{nk + i_{1}k}(\Theta(q))$ and
$\phi^{-m}(A_{0}) \in h^{i_{r}k}(\Theta(q))$ for all $m \geq 1$.
These have projections which are a distance at least $5$ apart, unless
$n+i_{1} = i_{r}$.
However, in this case $\phi$ is conjugate to
$$h^{nk} \phi_{i_{2}} \cdots \phi_{i_{r-1}} (\phi_{i_{r}} h^{nk}
\phi_{i_{1}} h^{-nk})$$
and $\phi_{i_{r}} h^{nk} \phi_{i_{1}} h^{-nk} \in
h^{i_{r}k}G(q)h^{-i_{r}k}$ which contradicts the minimality of $r$.
We may now proceed as in Theorem \ref{combothrm}. $\Box$

%%%%%%%%%%%%%%%%%%%%%%%%%%%%%%%%%%%%%%%%%%%%%%%
%%%%%%%%%%%%%%%%%%%%%%%%%%%%%%%%%%%%%%%%%%%%%%%
%%%%%%%%%%%%%%%%%%%%%%%%%%%%%%%%%%%%%%%%%%%%%%%
\section{Concluding remarks} \label{remarksect}
%%%%%%%%%%%%%%%%%%%%%%%%%%%%%%%%%%%%%%%%%%%%%%%
%%%%%%%%%%%%%%%%%%%%%%%%%%%%%%%%%%%%%%%%%%%%%%%
%%%%%%%%%%%%%%%%%%%%%%%%%%%%%%%%%%%%%%%%%%%%%%%

\subsection{} A more geometric description of Veech groups is the following.
Any $q \in {\cal A}_{2}{\cal T}(R)$ defines a holomorphic totally geodesic
embedding of the hyperbolic plane
$$f_{q}: {\mathbb H}^{2} \rightarrow {\mathbb H}_{q} \subset {\cal T}(R)$$
The stabilizer $\Stab_{\Mod(R)}({\mathbb H}_{q})$ acts on ${\mathbb
H}_{q}$ and this action can be conjugated back to ${\mathbb H}^{2}$, via
$f_{q}$, thus defining a homomorphism to $\PSL_{2}({\mathbb R}) \cong
\Aut({\mathbb H}^{2})$.
It turns out that $\Stab_{\Mod(R)}({\mathbb H}_{q}) = \Aff^{+}(q)$ and
with the appropriate identifications, the homomorphism to
$\PSL_{2}({\mathbb R})$ is given by $D$.
In the special case that ${\mathbb H}^{2}/\PSL(q)$ has finite area, this
quotient is called a Teichm\"uller curve and this immerses into the moduli
space.

This sheds some additional light on the analogy with Kleinian groups.
In particular, when the Fuchsian subgroups $G_{1},G_{2}$ lie in a fixed
Kleinian group $\Gamma$, these define totally geodesic surfaces in the
associated hyperbolic $3$-manifold ${\mathbb H}^{3}/\Gamma$.
If $h$ is taken to lie in $\Gamma$ (and is sufficiently complicated), then
the amalgamated product injects into $\Gamma$.
In this case, one can view the construction as truncating one of the cusps
of each of the totally geodesic surfaces and connecting the exposed
boundary components with an annulus whose co-core is determined by $h$.
Similarly, we can view our construction as truncating one of the cusps of
the Teichm\"uller curves and connect the exposed boundary components by
an annulus in the moduli space whose co-core is determined by $h$.
We remark that our theorem (as in the Kleinian group case) does not
require $\PSL(q)$ to have finite area.\\

To further the analogy with hyperbolic spaces, we recall that Thurston's
compactification of ${\cal T}(R)$ is obtained by adding $\prj\ML_{0}(R)$
at infinity to obtain
$$\overline{{\cal T}}(R) = {\cal T}(R) \cup \prj\ML_{0}(R) \cong
B^{6g-6}$$
where $B^{6g-6}$ is the closed ball of dimension $6g-6$ and
$\prj\ML_{0}(R)$ is identified with the boundary.
Moreover, $\prj\ML_{0}(R)$ has a natural piecewise projective structure
and the action of $\Mod(R)$ on ${\cal T}(R)$ and $\prj\ML_{0}(R)$ fit
together to give a well defined action on $\overline{\cal T}(R)$ which is
holomorphic on the interior and piecewise projective on the boundary.

There is a natural identification of $\prj \ML_{0}(q)$ with the boundary
at infinity $\partial_{\infty}{\mathbb H}^{2}$.
In this way, the inclusion of $\prj \ML_{0}(q)$ into $\prj \ML_{0}(R)$ can
be thought of as an extension $\partial_{\infty} f_{q}$ of $f_{q}$ to
infinity.
Indeed, the natural projective structure on ${\mathbb R}\prj^{1} =
\partial_{\infty}{\mathbb H}^{2} = \prj \ML_{0}(q)$ makes
$\partial_{\infty}f_{q}$ into a piecewise projective embedding,
equivariant with respect to the $\Aff^{+}(q)$ action.

Moreover, $\partial_{\infty} f_{q}$ sends the limit set $\Lambda(\PSL(q))
\subset \partial_{\infty} {\mathbb H}^{2} = \prj \ML_{0}(q)$
homeomorphically and $\Aff^{+}(q)$-equivariantly to the limit set
$\Lambda(\Aff^{+}(q)) \subset \prj \ML_{0}(R)$, as defined by McCarthy and
Papadopoulos \cite{MP}.

The map
$$\overline{f}_{q} = f_{q} \cup \partial_{\infty} f_{q}: \overline{\mathbb
H}^{2} = {\mathbb H}^{2} \cup \prj \ML_{0}(q) \rightarrow \overline{\cal
T}(R)$$
is continuous for every $p \in {\mathbb H}^{2}$ and almost every $p \in
\prj \ML_{0}(q)$ by a theorem of Masur \cite{Mas2}.
However, Masur's theorem implies that this is in general not continuous at
{\em every} point of $\prj \ML_{0}(q)$.

This now leads us to the following natural question, the analog of which
is true in the setting of Kleinian groups.

\begin{question}
Let $G \cong \pi_{1}(S_{2g}) \rightarrow \Mod(R_{g})$ be the embedding
given by Corollary \ref{sfcegrpcor}.
Consider $\partial_{\infty}(G)$ which can be canonically identified with
the circle at infinity of the universal cover $\widetilde{S}_{2g} \cong
{\mathbb H}^{2}$ of $S_{2g}$.
Does there exist a continuous $G$-equivariant map
$$\partial_{\infty}(G) \rightarrow \prj \ML_{0}(R)$$
\end{question}

It is not hard to see that if such a map exists, it must be unique.
In fact, the proposed map is already defined on a subset of
$\partial_{\infty}(G)$, namely, the limit points of the conjugates of the
Veech groups in the amalgamation.

\subsection{} It is known that many classes of lattices in 
Lie groups cannot inject into $\Mod(R)$. Indeed if the lattices are
superrigid the image of such a lattice in $\Mod(R)$ is necessarily
finite (see \cite{FM} and \cite{Y}).  From this it follows
(cf. Theorem 2 of \cite{Y}), that the only lattices that can admit a
faithful representation (or even an infinite representation) into
$\Mod(R)$ are lattices in $\SO(m,1)$, $m\geq 2$ or $\SU(q,1)$,
$q\geq1$.  Indeed, since solvable subgroups of $\Mod(R)$ are virtually
abelian \cite{BLM}, this observation also excludes non-cocompact
lattices of $\SU(q,1)$ for $q\geq 2$ from injecting. 

In light of these comments and the results of this paper, it seems
interesting to ask whether such injections can arise for (cocompact)
lattices in $\SO(m,1)$, $m\geq 3$ or $\SU(q,1)$, $q\geq 2$.  There are
simple obstructions to injecting certain of these lattices, or indeed
for any group. Namely if a finitely generated group $G$ admits an
injection into $\Mod(R)$ for some $R$, the $\vcd(G) \leq
\vcd(\Mod(R))$ (see \cite{Br}).  If the vcd of a group $G$ satisfies
the above inequality, then we call $G$ {\em admissible for} $R$ or
simply {\em admissible}.  It follows that for lattices in $\SO(m,1)$,
$m\geq 3$ or $\SU(q,1)$, $q\geq 2$ and for a fixed $R$, the values of
$m$ and $q$ above are bounded.

The vcd's of the groups $\Mod(R)=\Mod(R_{g,n})$ are well-known (see \cite{Ha}). If 
$R$ is not sporadic, the vcd of $\Mod(R)$ is $4g-5$ if $n=0$,
$4g-4+n$ if $n>0$ and $g\geq 1$ and $n-3$ if $n\geq 5$.

Motivated by this discussion we pose.

\begin{question}
Suppose $R=R_{g,n}$ is not sporadic.  Let $\Gamma$ be a lattice in 
$\SO(m,1)$, $m\geq 3$ or $\SU(q,1)$, $q\geq 2$, admissible for $R$. Does $\Gamma$ 
inject in $\Mod(R)$?
If such an injection exists does there exist 
a continuous $\Gamma$-equivariant map
$$\partial_{\infty}(\Gamma) \rightarrow \prj \ML_{0}(R)?$$
\end{question}

More generally, for a fixed $R$ and hence fixed vcd, one can
ask the above question for negatively curved groups in the sense of Gromov.  

Thus we pose the following question.

\begin{question}
Suppose $R=R_{g,n}$ is not sporadic.  Which 1-ended
negatively curved groups $G$, admissible for $R$, inject in $\Mod(R)$?
If such an injection exists does there exist 
a continuous $G$-equivariant map
$$\partial_{\infty}(G) \rightarrow \prj \ML_{0}(R)?$$
\end{question}

For the sporadic surfaces, the associated mapping class groups are
virtually free, hence no $1$-ended group embeds.

Work of a similar vein is done in \cite{FMo}, where convex cocompact subgroups
of $\Mod(R)$ are investigated.

%%%%%%%%%%%%%%%%%%%%%%%%%%%%%%%%%

\begin{tabular}{lll}
Department of Mathematics & \hspace*{1truein} & Department of Mathematics\\
Columbia University & \hspace*{1truein} &The University of Texas at Austin\\
2990 Broadway MC 4448 & \hspace*{1truein} &1 University Station C1200\\
New York, NY 10027-6902 & \hspace*{1truein} &Austin, TX 78712-0257\\
email: {\tt clein@math.columbia.edu} & \hspace*{1truein} &email: {\tt areid@math.utexas.edu}\\ \end{tabular}

%Phone: (512) 471-3153\\
%Phone: (212) 854-2431\\
\end{document}